\documentclass[12pt]{amsart}

\usepackage{amsmath, amssymb, color} 
\usepackage{amsthm}
\usepackage{multicol}
\usepackage{amscd}
\usepackage[colorlinks,allcolors=blue]{hyperref}
\usepackage[nameinlink,noabbrev,capitalise]{cleveref}
\usepackage{tikz}
\usepackage{graphicx}
\usepackage{standalone}
\usepackage{geometry}
\usepackage{algorithm}
\usepackage{algpseudocode}
\usepackage{comment}
\usepackage{mathrsfs}  

\numberwithin{equation}{section}
\newtheorem{theorem}{Theorem}[section]
\newtheorem{lemma}[theorem]{Lemma}
\newtheorem{proposition}[theorem]{Proposition}
\newtheorem{corollary}[theorem]{Corollary}

\newtheorem{question}[theorem]{Question}
\newtheorem*{thm}{Theorem}

\theoremstyle{definition}
\newtheorem{definition}[theorem]{Definition} 
\newtheorem{remark}[theorem]{Remark}
\newtheorem{example}[theorem]{Example}

\DeclareMathOperator{\conv}{conv}
\DeclareMathOperator{\indeg}{indeg}
\DeclareMathOperator{\sd}{sd}

\newcommand{\cD}{\mathcal{C}}

\usetikzlibrary{decorations.markings}
\usetikzlibrary {arrows.meta,bending,positioning}
\usepackage{subcaption}
\usepackage{customcommands}

\usetikzlibrary{arrows.meta}

\begin{document}

\title{Parking functions and chip-firing on hypergraphs}
\date{\today}

\author[T. Blanton]{Timothy Blanton}
\address[T. Blanton]{Department of Mathematics, 
University of California, 
Davis, CA, USA}
\email{trblanton@ucdavis.edu}

\author[A. Dochtermann]{Anton Dochtermann}
\address[A. Dochtermann]{Department
of Mathematics,
Texas State University, TX, USA}
\email{dochtermann@txstate.edu}

\author[I. Hong]{Isabelle Hong}
\address[I. Hong]{
Department of Mathematics, University of California, Los Angeles, CA, USA}
\email{isabellehong@ucla.edu}

\author[S. Oh]{SuHo Oh}
\address[S. Oh]{Department
of Mathematics,
Texas State University, TX, USA}
\email{s\_o79@txstate.edu}

\author[Z. Zhan]{Zhan Zhan}
\address[Z. Zhan]
{Department of Mathematics, University of Washington, WA, USA}
\email{zzhan4@uw.edu}

\begin{abstract}
For a connected graph $G$ with sink vertex $q$, a $G$-parking function is a vector of nonnegative integers whose entries are determined by cut-sets in $G$. Such objects also arise as the superstable configurations in the context of chip-firing.  The set of all $G$-parking functions have various algebraic and combinatorial properties; for instance they relate to evaluations of the Tutte polynomial and in particular are counted by spanning trees of $G$.
We extend these constructions to the setting of hypergraphs, where edges can have multiple vertices.  For a hypergraph $H$ with sink $q$, we define $H$-parking functions in terms of cuts in $H$ and prove that the maximal such sequences are characterized by certain acyclic orientations of $H$. We introduce a notion of a $q$-rooted spanning tree for $H$, and prove that the set of all such objects are counted by $H$-parking functions. 
We also show how $H$-parking functions can be recovered as the superstable configurations in a version of chip-firing on $H$, where chips have a choice of where to go when fired. We prove that one can recover such configurations via chip-firing on a family of digraphs associated to $H$.
\end{abstract}

\maketitle

\section{Introduction}
Suppose $G = (V,E)$ is a connected graph with specified sink vertex $q \in V$ and nonsink vertices $[n] = \{1, 2, \dots, n\}$. A $G$-parking function is a sequence $\vec{c} \in {\mathbb Z}_{\geq 0}^n$ of nonnegative integers whose entries are bounded above by parameters associated to vertex cuts in $G$.  Such objects were studied by Postnikov and Shapiro in \cite{PostnikovShapiro03}, and when $G= K_{n+1}$ is the complete graph they recover the classical parking functions that originated in queing theory and have since been the subject of extensive study \cite{Yan15}. In this case, such objects have an easy description: a sequence is a parking function if and only if its nondecreasing rearrangement is coordinatewise bounded by the vector $\vec{v} = (1,2, \dots, n)$. This naturally leads to a more general notion of a \emph{vector parking function} for other choices of $\vec{v}$.

$G$-parking functions can also be understood in terms of \emph{chip-firing} on the graph $G$, where the distribution rule is determined by the reduced Laplacian $L_G = L_G^q$. Here we consider configurations $\vec{c} \in \mathbb{Z}_{\geq 0}^{n}$ of chips on the nonsink vertices of $G$. If a nonsink vertex has at least as many chips as its degree, it can fire, passing one chip to each of its neighbors.  If no non-sink vertices can be fired, we say that the configuration is \emph{stable}. 
A configuration $\vec{c}$ is said to be \emph{superstable} if no nonempty set of non-sink vertices can fire, so that subtracting from $\vec{c}$ any nonempty collection of the columns of $L_G$ results in a vector with a least one negative coordinate. For undirected graphs (and, more generally, Eulerian digraphs), the set of superstable configurations can be seen to coincide with the collection of $G$-parking functions.

According to the definitions, to determine whether a given $\vec{c} \in \mathbb{Z}_{\geq 0}^{n}$ is a $G$-parking function, one must test a condition for all nonempty subsets $T \subset V \setminus q$.
However, Dhar's burning algorithm \cite{Dhar90} provides a polynomial-time way to check this condition.
 $G$-parking functions also have many pleasing combinatorial properties, and for instance the collection of all such sequences is in bijection with the set of spanning trees of $G$.
In \cite{ChebikinPylyavskyy05}, Chebikin and Pylyavskyy describe a family of such bijections between $G$-parking functions and spanning trees. They work in the more general setting of digraphs, where the relevant combinatorial objects are $q$-rooted spanning trees. For the case of undirected graphs, a particular specialization of their construction can be seen to recover Dhar's algorithm.

Varying the choice of sink vertex $q$ will in general lead to a different set of $G$-parking functions for a given graph $G$. However, some parameters are invariant. As mentioned above, if $G$ is connected then the set of $G$-parking functions is in bijection with the set of spanning trees of $G$. Furthermore, for a graph $G$ with a specified sink $q$, the \emph{maximal} $G$-parking functions are in bijection with acyclic orientations of $G$ with a unique source $q$ \cite{BensonChakrabartyTetali10}. From this it follows that all maximal $G$-parking functions have the same degree. In \cite{Merino97}, Merino showed that the \emph{degree sequence} of $G$-parking functions can be recovered as an evaluation of the Tutte polynomial of $G$ (and hence is independent of the choice of $q$). This then can be used \cite{Merino01} to establish a conjecture of Stanley regarding the $h$-vector of matroids for the special case of cographic matroids.

\subsection{Our contributions}
In this work we seek to extend the constructions and results described above to the setting of \emph{hypergraphs}, where edges can now have more than one `endpoint'.
For this we fix a hypergraph $H = (V,E)$ with a specified sink vertex $q \in V$ and nonsink vertices $[n] = \{1, 2, \dots, n\}$.
We let $B(H)$ denote its bipartite incidence graph, with vertex set ${\mathcal V} = (E,V)$ and edges $\{e,v\} \in {\mathcal E}$ whenever $v \in e$. Many of our constructions on hypergraphs will be phrased in terms of the graph $B(H)$.

For any nonempty $T \subset [n]$ and $v \in T$, we define the \emph{degree} of $v$ (with respect to $T$) to be the number of hyperedges that contain $v$, and which are not contained in $T$. Note that this notion generalizes the concept of a cut set for classical graphs. As in the graph setting, we then define an $H$-parking function to be a sequence $\vec{c} = (c_1, c_2, \dots, c_n)$ where for any $T \subset [n]$, we have some $i \in T$ where $c_i$ is strictly less than the degree of $i$ with respect to $T$ (see \cref{def:outdegree} and \cref{def:Hparking}). In \cref{prop:conversion} we show how $H$-parking functions on a hypergraph $H$ can be understood in terms of certain $G$-parking functions on the underlying incidence bipartite graph $B(H)$.

By definition, the set of $H$-parking functions is an order ideal (meaning that subtracting $1$ from an entry will result in another $H$-parking function, as long as that entry is nonnegative). To understand the \emph{maximal} $H$-parking functions, in \cref{defn:orientation} we define the notion of an orientation of a hypergraph $H$. This leads to a combinatorial characterization of maximal $H$-parking functions, analogous to the situation for graphs.

\begin{thm}[\cref{thm:maximalsup}]
There exists a bijection between the set of acyclic orientations of $H$ with unique source $q$ and the set of maximal $H$-parking functions on $H$.
\end{thm}

As a consequence, we see that for any hypergraph $H$, the degree sequence set of $H$-parking functions forms a pure $O$-sequence.
We next study $H$-parking functions for \emph{complete hypergraphs} $K_{n+1}^d$, by definition hypergraphs whose edge set consists of all $d$-subsets of $[n+1]$. It turns out that such sequences can be recovered as vector parking functions.

\begin{thm}[\cref{thm:ssuparking}]
The $H$-parking functions on $K_{n+1}^d$ can be described as ${\vec u}$-parking functions, for specified values of ${\vec u} = (u_1, u_2, \dots, u_n)$.
\end{thm}

Recall that in the classical case, there exists a bijection between the set of spanning trees of a connected graph $G$ and the set of $G$-parking functions. In a similar way, we wish to identify underlying combinatorial structures in a hypergraph $H$ that count $H$-parking functions. In \cref{def:tree} we define the notion of a \emph{$q$-rooted spanning tree} of a rooted hypergraph $H$, based on a notion of burning equivalence for spanning trees of the underlying bipartite graph $B(H)$.  We then have the following result.

\begin{thm}[\cref{thm:bijection}]
For any hypergraph $H$ with sink $q$, there exists a bijection between the set of $H$-parking functions of $H$ and the set of $q$-rooted spanning trees of $H$.
\end{thm}

We emphasize that our notion of a spanning tree of $H$ is dependent on the choice of sink vertex $q$. Indeed, the number of $H$-parking functions can change when one varies this choice. Our definition of spanning trees also relates to other constructions from the literature, and in \cref{prop:hypertree} we show how our spanning trees refine the notion of a \emph{hypertree} of a hypergraph introduced by K\'alma\'n and Postnikov in \cite{Kalman17}.

\subsection{Chip-firing}

We next investigate how $H$-parking functions can be understood in terms of a version of \emph{chip-firing} on hypergraphs. As usual, we consider a configuration $\vec{c} \in {\mathbb Z}_{\geq 0}^n$ on the nonsink vertices of $H$. As in the classical setting, when we fire a vertex $v$, we can think of it as sending a chip to each of its incident edges. For a usual graph (where edges have just two endpoints) there is a canonical place to send this chip. However, if an edge $e \in E(H)$ has cardinality at least 3, we now have a choice of where to send the chip. We make this notion precise in our definition of a \emph{firing choice}, see \cref{sec:hyper}.

To develop our notion of superstability we wish to fire subsets $T \subset V \setminus q$ of nonsink vertices, and here some subtleties arise. First, we insist that if $e \in E(H)$ is an edge such that $e \subset T$, then our firing choice restricted to $e$ must result in a permutation of the underlying set of chips (we say that such a choice is \emph{cancellative}, see \cref{def:setfire}). In addition, it is possible that a certain cancellative firing choice leads to a valid configuration, whereas another leads to a configuration with negative values. We say that $T$ is ready to fire if \emph{any} cancellative firing choice results in a nonnegative configuration, see \cref{def:setfire}.  We then define $\vec{c}$ to be superstable if no nonempty subset of vertices is ready to fire.
In \cref{lem:readytofire} we show that a configuration $\vec{c}$ on a hypergraph $H$ is superstable if and only if $\vec{c}$ is an $H$-parking function.

  As chip-firing on a hypergraph $H$ involves a firing choice at each step, there is no single Laplacian-type matrix that describes the theory. However, it turns out $H$ defines a \emph{family} of matrices that can be used to recover the superstable configurations of $H$. For this, recall that firing a vertex of $H$ requires making a choice of which incident vertex an edge sends a chip to. One way to make this choice is to fix a cyclic order on each edge. The resulting \emph{cycling} ${\mathcal O}$ defines an Eulerian digraph $D_{\mathcal O}(H)$ (see \cref{sec:digraphs}) for which the digraph Laplacian $L(D_{\mathcal O}(H))$ defines a chip-firing rule whose superstable configurations coincide with the $G$-parking functions for $D_{\mathcal O}(H)$. Our next result says that the superstable configurations for $H$ are exactly the superstable configurations obtained this way.

\begin{thm}[\cref{thm:hypersuperstable}]
    A configuration $\vec{c}$ is superstable for $H$ if and only if it is superstable for $D_{\cD}(H)$ for some choice of vertex-induced cycling $\cD$.
\end{thm}

It remains an open question to decide if the directed Laplacians defined by the various $D_{\cD}(H)$ can be used to understand the chip-firing rule for the hypergraph $H$.

\subsection{Related work}
In recent years, the theory of chip-firing on graphs has been generalized to other settings in a number of ways. Much of this work relies on the fact that chip-firing moves on a graph $G$ are described by the graph Laplacian $L_G$. Recall that $L_G = i i^T$, where $i$ is the (reduced) signed incidence matrix, so that $L_G$ is an $n \times n$ matrix whose $i$th row corresponds to firing the vertex $v_i$. Using these ideas, Duval, Klivans, and Martin \cite{DuvKliMar} introduced a notion of higher dimensional chip-firing, where chip-firing moves are defined by the Laplacian associated to the ridge-facet incidences of a pure simplicial complex.  Here chips are placed on the ridges (codimension one faces) and passed along facets. This can be generalized to other cell complexes, and for instance in \cite{FelKli} Felzenszwalb and Klivans study the dynamical properties of flow firing on two dimensional complexes. Although one can naturally think of a hypergraph as a simplicial complex (eg by using the hyperedges to define facets), far as we know our work has no connection to this theory.

In another direction, Backman \cite{Backman} has introduced a notion of an oriented incidence matrix for a $d$-uniform hypergraph $H$. Again this is an `edge by vertex' matrix, but now the columns take entries among the $d$th roots of unity. This generalizes the case $d=2$, where each column has a single $1$ and a single $-1$ entry.  This matrix then leads to a Laplacian $L_H$ that one can use to define a chip-firing rule on $H$, although again we do not know of any connection to our work.

Finally we wish to mention work of Dong from \cite{Dong}, where the notion of a \emph{$B$-parking function} on a bipartite graph $B$ is defined. As in our work, the combinatorial objects that count $B$-parking functions are generalizations of spanning trees of a graph $G$, expressed in terms of its bipartite incident graph $B(G)$.  As noted in \cite{Dong}, if $G$ is a graph with vertex set $V$ and root $q$, there is a bijection between the set of spanning trees of $G$ and the set of uniquely restricted matchings of size $|V|-1$ in $B(G) \backslash q$. Here a matching is \emph{uniquely restricted} if it is the only perfect matching in the subgraph of $G$ induced by the vertices of $M$.
Motivated by this, Dong defines a notion of a \emph{$B$-parking function} for an arbitrary bipartite graph $B$, and establishes a bijection between the set of uniquely restricted matchings of $B$ and the set of $B$-parking functions. 
Although our two notions agree when $B = B(G)$ is the bipartite incidence graph of a (classical $2$-regular) graph $G$, one can see that in general the two theories diverge. 

\subsection{Organization}
The rest of the paper is organized as follows. In \cref{sec:digraphbackground} we review the necessary background from the theory of parking functions and chip-firing on graphs. Here we work in the more general context of digraphs, and recall a bijection between digraphs and $q$-rooted trees. In \cref{sec:hyper} we turn to hypergraphs, and provide the definition of an $H$-parking function for a hypergraph $H$. In \cref{sec:maximal} we define the notion of an orientation of a hypergraph (relative to a choice of sink vertex $q$), and use these to provide a characterization of maximal $H$-parking functions. In \cref{sec:complete} we consider \emph{complete} hypergraphs, and show how $H$-parking in this context can be recovered as vector parking functions. In \cref{sec:bijections} we define our notion of ($q$-rooted) spanning trees of a hypergraph $H$ and prove that such objects are in bijection with the set of $H$-parking functions.
In \cref{sec:stars} we consider the special case of star hypergraphs, where the sink vertex $q$ appears in every hyperedge of $H$. In this case we show that the chip-firing theory reduces to chip-firing on a certain digraph.

In \cref{sec:chipfiring} we consider $H$-parking functions in the context of a theory of chip-firing on the hypergraph $H$. Here we show that the notion of a superstable configuration coincides with that of an $H$-parking function. In \cref{sec:digraphs} we show how a choice of cyclings of the edges of $H$ lead to a collection of digraphs that can be used to recover the $H$-parking functions. In \cref{sec:commalg} we briefly discuss how $H$-parking functions can also be understood in terms of generators and resolutions of certain $H$-parking function ideals (as studied in \cite{AlmDocSmi}). Finally, in \cref{sec:further} we discuss some open questions and potential further directions of study.

\section{Parking functions and chip-firing on digraphs}\label{sec:digraphbackground}

We begin by recalling some notions from the theory of parking functions and chip-firing. We work in the setting of directed graphs (digraphs), since we will need these more general constructions in our study of hypergraphs. We discuss the basic concepts here and refer to \cite{Holroyd08}, \cite{Klivans19}, and \cite{CorryPerkinson18} (and the references therein) for details.

For us a \emph{directed graph} (or \emph{digraph}) $G = (V,E)$ consists of a vertex set $V$ and a multiset of directed edges $E \subset V \times V$. If $e = (v,w) \in E$ we say that $e$ is an edge from $v$ to $w$. In the context of digraphs we define the outdegree and indegree of a vertex $v \in V$ as 
\[\outdeg(v) = |\{e \in E: e = (v,w) \text{ for some $w \in V$}\}, \]
\[\indeg(v) =  |\{e \in E: e = (u,v) \text{ for some $u \in V$}\}. \]
\noindent
We say that $G$ is \emph{symmetric} if $(i,j) \in E$ if and only if $(j,i) \in E$. Symmetric graphs recover the notion of an \emph{undirected} graph, where we often simply write $\{i,j\} \in E$.

For our chip-firing setup, we fix a sink vertex $q \in V$ and let $[n] = \{1, 2, \dots, n\}$ denote the nonsink vertices.  We let $L_G = L_G^q$ denote the \emph{directed (reduced) Laplacian}, the $n$-by-$n$ matrix whose diagonal entries are given by the outdegree of the corresponding vertex, and whose off diagonal entries are given by $\ell_{i,j} = -| \{\text{edges from $i$ to $j$}\}|$. A \emph{configuration of chips} $\vec{c} \in {\mathbb Z}^n_{\geq 0}$ is a vector of nonnonegative integers indexed by the nonsink verties of $G$. Given a configuration $\vec{c}$, a vertex $i \in n$ is \emph{ready to fire} if $\outdeg(v) \geq c_i$, in which case we can fire $i$ to obtain the configuration
\[\vec{d} = \vec{c} - L_G \vec{e}_i,\]
where $\vec{e}_i$ is the corresponding standard basis vector. Hence $i$ is ready to fire if subtracting the $i$th column of $L_G$ from $\vec{c}$ results in a vector with nonnegative entries. Similarly, a nonempty subset $S \subset [n]$ is ready to fire if subtracting the corresponding set of columns of $L_G$ from $\vec{c}$ has nonnegative entries. A configuration is \emph{superstable} if no nonempty set $S \subset [n]$ is ready to fire.

In \cite{PostnikovShapiro03}, Postnikov and Shapiro introduced the concept of a $G$-parking function associated to a digraph $G$.  Again we fix a sink vertex $q \in V$ and let $[n] = \{1,2, \dots, n\}$ denote the non-sink vertices. A \emph{$G$-parking function} is then a sequence $(b_1, \dots, b_n)$ of non-negative integers with the property that for each nonempty subset $S \subset [n]$ of nonsink vertices of $G$, there exists a vertex $j \in S$ such that the number of edges from $j$ to vertices outside of $S$ is greater than $b_j$. If $G$ is symmetric, this recovers the notion of a superstable configuration on undirected graphs. For the undirected complete graph $G = K_{n+1}$ (and any choice of sink), these sequences coincide with the classical parking functions.

A digraph $G = (V,E)$ is \emph{Eulerian} if it is strongly connected and $\indeg(v) = \outdeg(v)$ for all $v \in V$. For example, any symmetric digraph is Eulerian.
Gabrielov \cite{Gab} has shown that if $G$ is Eulerian, then the set of superstable configurations of $G$ coincides with the set of $G$-parking functions. 
More generally, these two notions coincide whenever the matrix $L$ governing the chip-firing rule is an $M$-matrix with nonnegative row sums (see \cite{GuzmanKlivans15} for definitions and details). In this case, one also has that each equivalence class ${\mathbb Z}^n/\im (L_G)$ contains exactly one superstable configuration (with nonnegative entries, by definition).

\subsection{Bijections between spanning trees and parking functions}
\label{sec:burning}

In \cite{ChebikinPylyavskyy05}, Chebikin and Pylyavskyy describe a family of bijections between the set of $G$-parking functions and the set of ($q$-rooted) spanning trees of a digraph $G$.  To recall this, suppose $G$ has sink $q$ and nonsink vertices $[n]$. A \emph{subtree of $G$ rooted at $q$} is a subgraph $T$ containing $q$ with the property that for every vertex $v \in V(T)$, there exists a unique (directed) path in $T$ from $v$ to $q$. A subtree is \emph{spanning} if it contains all vertices of $G$. Given a $q$-rooted subtree of $G$, we will sometimes omit $q$ and simply call it a rooted subtree.

For each rooted subtree $T$, we let $\pi(T)$ be a total order on the vertices of $T$, and write $i <_{\pi(T)} j$ to denote that $i$ is smaller than $j$ in this order. A set of tree orders $\pi(T)$, for all rooted trees $T$ of $G$, is \emph{proper} if the following conditions hold, for all $T$:
\begin{enumerate}
 \item if $(j,i)$ is an edge of $T$, then $i <_{\pi(T)} j$;
 \item if $t$ is a subtree of $T$ rooted at $q$, then the order $\pi(t)$ is consistent with $\pi(T)$; in other words, $i<_{\pi(t)} j$ if and only if $i <_{\pi(T)} j$ for all $i,j \in t$.
\end{enumerate}

For any rooted tree $T$ and a vertex $j$ of $G$, the ordering $\pi(T)$ induces an ordering on the edges directed from $j$ to vertices of $T$: we say $(j,i) < (j,i')$ whenever $i <_{\pi(T)} i'$. As shown in \cite{ChebikinPylyavskyy05}, a proper set of tree orders gives rise to a bijection between the rooted spanning trees of $G$ and the $G$-parking functions: for a rooted spanning tree $T$ and a vertex $j \in [n]$, let $e_j$ be the edge of $T$ going out of $j$. To obtain our $G$-parking function $(b_1,\dots,b_n)$, we define $b_j$ to be the number of edges $e$ going out of $j$ such that $e <_{\pi(T)} e_j$.

\begin{theorem}\cite[Theorem 2.1]{ChebikinPylyavskyy05}\label{thm:CPtheorem}
The map described above is a bijection between the set of $q$-rooted spanning trees of $G$ and the set of $G$-parking functions of $G$.
\end{theorem}

The inverse map is also explicitly described in \cite{ChebikinPylyavskyy05}. If $(b_1, \dots, b_n)$ is a $G$-parking function, one constructs the corresponding spanning tree one edge at a time. We start with the tree $t_0$ consisting of the vertex $q$ and set $p_0 = 0$. For $m \geq 1$, we construct a $q$-rooted tree $t_m$ as follows. Let $U_m$ denote the set of vertices not in $t_{m-1}$, and let $V_m$ be the set of vertices $j \in U_m$ such that \[|\{\text{edges from $j$ to $t_{m-1}$}\}| \geq b_j +1.\]
For each $j \in V_m$, let $e_j$ be the edge from $j$ to $t_{m-1}$ such that exactly $b_j$ edges $e$ from $j$ to $t_{m-1}$ satisfy $e < e_j$ in the order $\pi(t_{m-1})$. Let $t$ be the tree obtained by adjoining each vertex $j \in V_m$ to $t_{m-1}$ using the edge $e_j$. Let $p_m$ denote the smallest vertex of $V_m$ in the order $\pi(t)$, and set $t_m$ to be the tree obtained by connecting $p_m$ to $t_{m-1}$ with the edge $e_{p_m}$. We continue this process until $m = n$, and set $T = t_n$.

The bijection described in \cref{thm:CPtheorem} depends on a choice of (a set of) tree orders. For our purposes, we focus on a particular order called the \emph{breadth-first search order} in \cite{ChebikinPylyavskyy05}. First fix a total order $\beta$ on the vertex set of $G$. For a $q$-rooted tree $T$ and vertex $i \in T$, we define the height $h_T(i)$ of $i$ in $T$ to be the number of edges in the unique path from $i$ to the root $q$. We set $i <_{\pi(T)} j$ if $h_T(i) < h_T(j)$ or if $h_T(i) = h_T(j)$ and $i <_\beta j$. Then, as shown \cite{ChebikinPylyavskyy05}, the set $\pi(T)$ is a proper set of tree orders. In what follows, the bijection described in \cref{thm:CPtheorem} coming from the breath-first search order will be called the \emph{Breadth-first Chebikin Pylyavskyy (BCP) algorithm}. In Section \ref{sec:bijections}, we will adapt this algorithm to the setting of hypergraphs.

\section{Hypergraphs and \texorpdfstring{$H$}--parking functions} \label{sec:hyper}
\begin{figure}[h]
\centering
\begin{multicols}{2}
    \begin{tikzpicture}[scale=1]
    \definecolor{color1}{RGB}{200,0,0}
    \definecolor{color2}{RGB}{0,0,200}
    \definecolor{color3}{RGB}{25,161,21}
    \coordinate [label=left:\textcolor{black}{$4$}] (A) at (0,0);
    \coordinate [label=right:\textcolor{black}{$1$}] (B) at (0,2);
    \coordinate [label=left:\textcolor{black}{$2$}] (C) at (2,2);
    \coordinate [label=right:\textcolor{black}{$3$}] (D) at (2,0);
    \coordinate (E) at (1.3,1.3);
    \coordinate (F) at (0.7,1.3);
    \coordinate (G) at (0.7,0.7);
    \filldraw [black] (A) circle [radius=1.5pt]
               (B) circle [radius=1.5pt]
               (C) circle [radius=1.5pt]
               (D) circle [radius=1.5pt];
    \draw[color1, rotate=-45] (E) circle (60pt and 40pt);
    \draw[color3, rotate=45] (F) circle (60pt and 40pt);
      \draw[color2, rotate=-45] (G) circle (60pt and 40pt);
\end{tikzpicture}
        \label{fig:examplehyper}   
    \columnbreak
    \vspace*{.1 in}
                \begin{tikzpicture}[scale = 1]
                \definecolor{color1}{RGB}{200,0,0}
                \definecolor{color2}{RGB}{25,161,21}
                \definecolor{color3}{RGB}{0,0,200}
                \coordinate (D) at (2,0);
                \coordinate [label=right:\textcolor{black}{$3$}] (C) at (2,1);
                \coordinate [label=right:\textcolor{black}{$2$}] (B) at (2,2);
                \coordinate [label=right:\textcolor{black}{$1$}] (A) at (2,3);
                \coordinate [label=right:\textcolor{black}{$4$}] (D') at (2,0);
                \coordinate [label=left:\textcolor{black}{$134$}] (E) at (0,0.5);
                \coordinate [label=left:\textcolor{black}{$124$}] (F) at (0,1.5);
                \coordinate [label=left:\textcolor{black}{$123$}] (G) at (0,2.5);
                \draw [color1] (G) to (C)
                        (G) to (B) (G) to (A);
                \draw [color2] (F) to (B)
                        (F) to (D) (F) to (A);
                \draw [color3] (E) to (A)
                        (E) to (C) (E) to (D);
                \filldraw [black] (A) circle [radius=1.5pt]
                           (B) circle [radius=1.5pt]
                           (C) circle [radius=1.5pt]
                           (D) circle [radius=1.5pt]
                           (E) circle [radius=1.5pt]
                           (F) circle [radius=1.5pt]
                           (G) circle [radius=1.5pt];
            \end{tikzpicture}
        \label{fig:example of bipartite graph}
\end{multicols}
\caption{A hypergraph $H$ along with its bipartite representation $B(H)$.}
\label{fig:bipartite}
\end{figure}
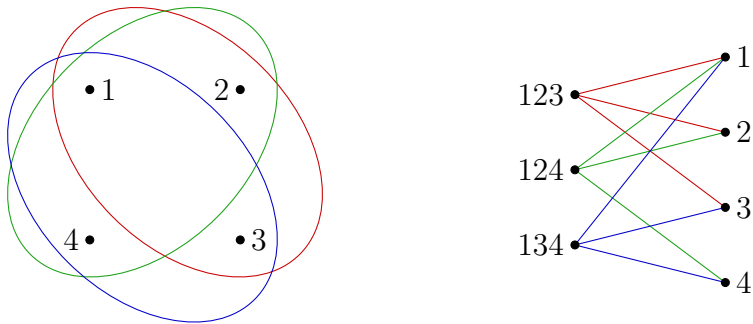

A \emph{hypergraph} $H = (V,E)$ consists of a vertex set $V = V(H)$ and a set $E = E(H)$ of subsets of $V$ called hyperedges. Elements of $E$ will be called \emph{edges} if the context is clear. We can represent a hypergraph $H$ in terms of its bipartite incidence graph $B(H) = (\cV, \cE)$, with vertex set $\cV = E \sqcup V$ and with edges $(e,v) \in \cE$ whenever $v \in e$. We refer to \cref{fig:bipartite} for a running example that we will use throughout the paper. If all elements of $E$ have the same cardinality $d$, we say that $H$ is \emph{$d$-regular}. 
The \emph{complete hypergraph} $K_n^d$ is the $d$-regular hypergraph on vertex set $[n] = \{1,2, \dots, n\}$ whose edge set consists of all $d$-subsets of $[n]$.

The \emph{degree} of a vertex $v \in V(H)$ is given by
\[\deg(v) = |\{e \in E(H): v \in e\}|.\]
Note that $\deg(v)$ agrees with the degree of $v$ in the (usual) graph $B(H)$.

\setlength\columnsep{-100pt}

\begin{remark}
Throughout the paper, we will assume that our hypergraphs $H$ are \emph{connected}, i.e., given any two vertices $v_1$ and $v_2$, there exists a sequence of edges $(e_1, e_2, \dots, e_m)$ such that $v_1 \in e_1$, $v_2 \in e_m$ and $e_i \cap e_{i+1} \neq \emptyset$ for all $i = 1, \dots m-1$. From this it follows that the graph $B(H)$ is also connected.
\end{remark}
    
We will also need a notion of degree of a vertex relative to a subset that contains it.

\begin{definition} \label{def:outdegree} Suppose $H = (V,E)$ is a hypergraph, and let $T \subset V$ be a nonempty set of vertices with $v \in T$. The \emph{$T$-degree of $v$} is 
\[\deg^H_T(v) = |\{e \in E: \text{$v \in e $ and $e \not \subset T$}\}|.\]
\end{definition}

If the underlying hypergraph is clear, we will write $\deg_T(v)$.
Also note that if $T = \{v\}$ then we have $\deg_{\{v\}}(v) = \deg(v)$. With this we can provide one of the main definitions of the paper.

\begin{definition}\label{def:Hparking}
Suppose $H = (V,E)$ is a hypergraph with sink vertex $q \in V$ and nonsink vertices $[n]$. An \emph{$H$-parking function} is a sequence $\vec{c} \in {\mathbb Z}^n_{\geq 0}$ such that for all $T \subset [n]$, there exists an $i \in T$ such that $c_i < \deg_T(i)$.
\end{definition}

With the connection to chip-firing in mind, we will often refer to a vector $\vec{c} \in {\mathbb Z}^n_{\geq 0}$ as a \emph{configuration of chips}.  If a nonempty subset $T \subset [n]$ has the property that $\deg_T(i) \leq c_i$ for all $i \in T$, we will say that $T$ is \emph{bounded} by $\vec{c}$. Hence a configuration $\vec{c}$ is an $H$-parking function if no nonempty subset of $[n]$ is bounded by $\vec{c}$.

It turns out that $H$-parking functions on a hypergraph $H$ can be understood in terms of its underlying bipartite representation $B(H)$. Recall that the vertex set of $B(H)$ is given by ${\mathcal V} = E \sqcup V$. If $q \in V$ is the chosen sink vertex for $H$, then by convention we take $q$ to be the sink vertex for $B(H)$. Also, if $\vec{c}$ is a configuration on the nonsink vertices of $H$, we let $\vec{c}_B$ denote the configuration on the nonsink vertices of $B(H)$ given by placing zeros on each vertex $w \in E$. With this we have the following observation. 

\begin{proposition}
\label{prop:conversion}
Suppose $\vec{c}$ is a configuration on the nonsink vertices of a hypergraph $H$. Then $\vec{c}$ is an $H$-parking function for $H$ if and only if $\vec{c}_B$ is a $G$-parking function for $B(H)$.
\end{proposition}

\begin{proof}
Let $\vec{c}$ and $\vec{c}_B$ be as above, and suppose $T \subset V \setminus q$ is a subset of nonsink vertices of $H$. First we claim that $T$ is bounded by $\vec{c}$ if and only if 
\[T_B := \{w \in E \mid w \subseteq T\} \sqcup T \]
\noindent
is bounded by $\vec{c}_B$. To see this, let $v \in T_B$. Note that if $v \in E$ then the degree condition is trivially satisfied since $\deg_{T_B}^{B(H)}(v) = 0$. On the other hand, if $v \in T$, note that $\deg_T^H(v) = \deg_{T_B}^{B(H)}(v)$. The claim follows. This also shows that if $\vec{c}_B$ is an $G$-parking function for $B(H)$ then $\vec{c}$ is a $H$-parking function for $H$.

Now suppose $S' \sqcup S$ is bounded by $\vec{c}_B$, where $S' \subseteq E$ and $S \subset V \setminus q$. Note that since $(\vec{c}_B)_w = 0$ for all $w \in S'$, we must have $w \subset S$. Hence, it follows that $S' \sqcup S \subset S_B$, where $S_B$ is defined as above. Since $S' \sqcup S$ was bounded  $\vec{c}_B$, so is $S_B$. To see this, let $w \in S_B$. Note that if $w \in V$ then the degree can only be decreased, and if $w \in E$ then since $w \subset S'$ we again have $\deg_{S_B}(w) = 0$.
The claim then follows. From the previous paragraph, we conclude that $S$ is bounded by $\vec{c}$ (recall that $S \subset V \setminus q$).  This shows that if $\vec{c}$ is an $H$-parking function for $H$, then $\vec{c}_B$ is a $G$-parking function for $B(H)$, which completes the proof.
\end{proof}

\begin{remark}
Suppose $G$ is a graph on vertex set $V$ with sink vertex $q \in V$. Let $\sd(G)$ denote the subdivision of $G$, by definition the graph obtained from $G$ by replacing each edge $\{u,v\} \in E$ with the edges $\{u, w\}$ and $\{w, v\}$ (and adding the vertex vertex $w$). We see that in this case we have $\sd(G) = B(G)$, and hence from \cref{prop:conversion} we have that the set of $G$-parking functions for $G$ coincide with the set of $G$-parking functions of $\sd(G)$ that have $0$ entries on all the newly added vertices.
\end{remark}

\begin{remark}
If $H$ is any hypergraph, note that from \cref{prop:conversion} we can use Dhar's burning algorithm on $B(H)$ to check whether a given configuration on $H$ is an $H$-parking function.
\end{remark}

\subsection{Orientations and maximal parking functions}\label{sec:maximal}

For a usual graph $G$ with sink $q$, there is a simple bijection between the set of maximal $G$-parking functions and the set of \emph{acyclic orientations} of $G$ with a unique source $q$, as described in \cite{BensonChakrabartyTetali10}. Here we extend this to the setting of hypergraphs. We begin with a definition.

\begin{definition}\label{defn:orientation}
    Suppose $H$ is a hypergraph with sink $q$, and with bipartite incidence graph $B(H)$. An \emph{orientation} of an edge $e \in E(H)$ is a choice of vertex $v_e \in e$. An \emph{orientation of $H$} is a choice of orientation for each edge, which we denote ${\mathcal O} = \{v_e\}_{e \in E(H)}$.  
\end{definition}

    A orientation ${\mathcal O}$ of $H$ defines an orientation of the bipartite graph $B(H) = E(H) \sqcup V(H)$ by prescribing the directed edge $(v,e)$ if $v = v_e$ is the orientation of $e$, and $(e,v)$ otherwise.
    
\begin{definition}
   An orientation ${\mathcal O}$ of a hypergraph $H$ is \emph{acyclic with unique source $q$} if the underlying orientation of $B(H)$ has this property.
\end{definition}

\begin{remark}
We note that our definition of a hypergraphic orientation is a special case of a more general notion introduced in \cite{BenBer}, where an orientation of an edge $e \in E(H)$ is given by an ordered partition $(a,b)$, where $a$ and $b$ are nonempty subsets of $e$ satisfying $e = a \cup b$ and $a \cap b = \emptyset$. In our setup we always take $|a| = 1$. In \cite{BenBerMac} the authors show how (acyclic) orientations can be used to count the faces of \emph{hypergraph polytopes}. These objects also make an appearance in our study, see \cref{sec:genperm}.
\end{remark}

\begin{figure}[h]
\centering
\begin{multicols}{3}
    \begin{tikzpicture}[decoration={markings, mark= at position 0.65 with {\arrow[line width=0.4mm]{stealth}}}, scale = 1]
        \definecolor{color1}{RGB}{200,0,0}
        \definecolor{color2}{RGB}{0,0,200}
        \coordinate (D) at (2,0);
        \coordinate [label=right:\textcolor{color1}{$0$}] (C) at (2,1);
        \coordinate [label=right:\textcolor{color1}{$1$}] (B) at (2,2);
        \coordinate [label=right:\textcolor{color1}{$2$}] (A) at (2,3);
        \coordinate [label=left:\textcolor{black}{}] (E) at (0,0.5);
        \coordinate [label=left:\textcolor{black}{}] (F) at (0,1.5);
        \coordinate [label=left:\textcolor{black}{}] (G) at (0,2.5);
        \filldraw [black] (A) circle [radius=1.5pt]
                   (B) circle [radius=1.5pt]
                   (C) circle [radius=1.5pt]
                   (D) circle [radius=1.5pt]
                   (E) circle [radius=1.5pt]
                   (F) circle [radius=1.5pt]
                   (G) circle [radius=1.5pt];
        \draw [color2] (D) circle [radius=5pt];
        \draw[color=black, postaction={decorate}] (G) to (A);
        \draw[color=black, postaction={decorate}] (G) to (B);
        \draw[color=black, postaction={decorate}] (C) to (G);
        \draw[color=black, postaction={decorate}] (F) to (A);
        \draw[color=black, postaction={decorate}] (F) to (B);
        \draw[color=black, postaction={decorate}] (D) to (F);
        \draw[color=black, postaction={decorate}] (E) to (A);
        \draw[color=black, postaction={decorate}] (E) to (C);
        \draw[color=black, postaction={decorate}] (D) to (E);         
    \end{tikzpicture}
    
    \columnbreak
    \begin{tikzpicture}[decoration={markings, mark= at position 0.65 with {\arrow[line width=0.4mm]{stealth}}}, scale = 1]
    \definecolor{color1}{RGB}{200,0,0}
    \definecolor{color2}{RGB}{0,0,200}
    \coordinate (D) at (2,0);
    \coordinate [label=right:\textcolor{color1}{$1$}] (C) at (2,1);
    \coordinate [label=right:\textcolor{color1}{$0$}] (B) at (2,2);
    \coordinate [label=right:\textcolor{color1}{$2$}] (A) at (2,3);
    \coordinate [label=left:\textcolor{black}{}] (E) at (0,0.5);
    \coordinate [label=left:\textcolor{black}{}] (F) at (0,1.5);
    \coordinate [label=left:\textcolor{black}{}] (G) at (0,2.5);
    \filldraw [black] (A) circle [radius=1.5pt]
               (B) circle [radius=1.5pt]
               (C) circle [radius=1.5pt]
               (D) circle [radius=1.5pt]
               (E) circle [radius=1.5pt]
               (F) circle [radius=1.5pt]
               (G) circle [radius=1.5pt];
    \draw [color2] (D) circle [radius=5pt];
    \draw[color=black, postaction={decorate}] (G) to (A);
    \draw[color=black, postaction={decorate}] (B) to (G);
    \draw[color=black, postaction={decorate}] (G) to (C);
    \draw[color=black, postaction={decorate}] (F) to (A);
    \draw[color=black, postaction={decorate}] (F) to (B);
    \draw[color=black, postaction={decorate}] (D) to (F);
    \draw[color=black, postaction={decorate}] (E) to (A);
    \draw[color=black, postaction={decorate}] (E) to (C);
    \draw[color=black, postaction={decorate}] (D) to (E);

\end{tikzpicture}
    
    \columnbreak
    \begin{tikzpicture}[decoration={markings, mark= at position 0.65 with {\arrow[line width=0.4mm]{stealth}}}, scale = 1]
    \definecolor{color1}{RGB}{200,0,0}
    \definecolor{color2}{RGB}{0,0,200}
    \coordinate (D) at (2,0);
    \coordinate [label=right:\textcolor{color1}{$1$}] (C) at (2,1);
    \coordinate [label=right:\textcolor{color1}{$1$}] (B) at (2,2);
    \coordinate [label=right:\textcolor{color1}{1}] (A) at (2,3);
    \coordinate [label=left:\textcolor{black}{}] (E) at (0,0.5);
    \coordinate [label=left:\textcolor{black}{}] (F) at (0,1.5);
    \coordinate [label=left:\textcolor{black}{}] (G) at (0,2.5);
    \filldraw [black] (A) circle [radius=1.5pt]
               (B) circle [radius=1.5pt]
               (C) circle [radius=1.5pt]
               (D) circle [radius=1.5pt]
               (E) circle [radius=1.5pt]
               (F) circle [radius=1.5pt]
               (G) circle [radius=1.5pt];
    \draw [color2] (D) circle [radius=5pt];
    \draw[color=black, postaction={decorate}] (A) to (G);
    \draw[color=black, postaction={decorate}] (G) to (B);
    \draw[color=black, postaction={decorate}] (G) to (C);
    \draw[color=black, postaction={decorate}] (F) to (A);
    \draw[color=black, postaction={decorate}] (F) to (B);
    \draw[color=black, postaction={decorate}] (D) to (F);
    \draw[color=black, postaction={decorate}] (E) to (A);
    \draw[color=black, postaction={decorate}] (E) to (C);
    \draw[color=black, postaction={decorate}] (D) to (E);  
\end{tikzpicture}
\end{multicols}
\caption{Acyclic orientations with corresponding maximal superstable configurations.}
\label{fig:acyclic-maximalsup}   
\end{figure}
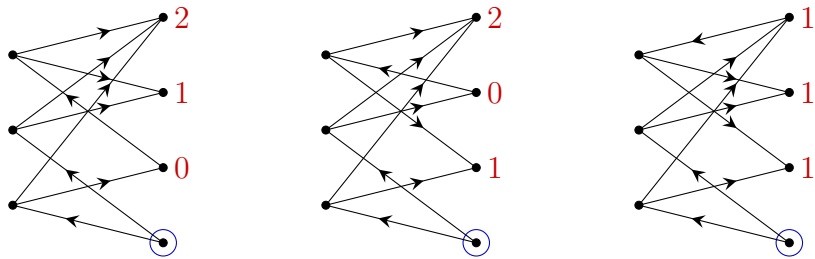

An orientation ${\mathcal O}$ of a hypergraph $H$ gives rise to a configuration $\vec{c} = \vec{c}({\mathcal O})$ on the nonsink vertices $[n] = \{1,2, \dots, n\}$ via the induced orientation on $B(H)$: For each $i \in [n]$ we let $\vec{c}_i = \indeg(i) -1$. See \cref{fig:acyclic-maximalsup} for an example.

\begin{theorem}\label{thm:maximalsup}
    Suppose $H$ is a hypergraph with sink $q$. Then the assignment described above defines a bijection between the set of acyclic orientations of $H$ with unique source $q$ and the set of maximal $H$-parking functions on $H$.
\end{theorem}

\begin{proof}
Recall from \cref{prop:conversion} that a vector $\vec{c}$  is a $H$-parking function if and only  $\vec{c}_B$ is a $G$-parking function for $B(H)$. From \cite{BensonChakrabartyTetali10} we know that the maximal $G$-parking functions on the graph $B(H)$ correspond to acyclic orientations of $B(H)$ with unique source $q$.  Here, for an orientation ${\mathcal O}$ of $B(H)$ we obtain a $G$-parking function $\vec{d}$ on $B(H)$ by taking ${\vec d}_i = \indeg_{\mathcal O}(i) - 1$.

This and \cref{prop:conversion} imply that among the maximal $G$-parking functions on $B(H)$,  the ones that correspond to $H$-parking functions on $H$ are those that come from orientations with indegree $1$ on every vertex $v \in E(H)$ (since these will result in 0 chips placed on every such $v$). Note that such orientations on $B(H)$ are precisely the acyclic orientations of $H$ with unique source $q$. The result follows.
\end{proof}

From the above result it follows that the degree of a maximal $H$-parking function equals the number of edges oriented rightwards minus the number of non-sink vertices of $H$, so it is independent of the acyclic orientation we choose.

\begin{corollary}
\label{cor:pure}
    For a hypergraph $H$ with chosen sink $q$, all maximal $H$-parking functions have the same degree.
\end{corollary}

From \cref{cor:pure} we have that the degree sequence of superstable configurations for any hypergraph $H$ is a \emph{pure $O$-sequence}.

\subsection{Complete hypergraphs and vector parking functions} \label{sec:complete}

We next consider $H$-parking functions for the case of $H = K_n^d$, the \emph{complete $d$-hypergraph on $n$ vertices} consisting of all $d$-subsets of $[n]$. In this section we show that the set of $H$-parking functions for such hypergraphs can be described as vector parking functions as studied by Yan \cite{Yan15}.

To recall this notion, for a sequence $(x_1, x_2, \dots, x_n)$ of real numbers, we let $x(1) \leq x(2) \leq \cdots \leq x(n)$ denote its rearrangement into a nondecreasing order.  Now we fix a nondecreasing vector of nonnegative integers ${\vec u} = (u_1, u_2, \dots, u_n)$. A vector $\vec{c} \in {\mathbb Z}^n_{\geq 0}$ is a $\emph{\bf u}$-parking function if its rearrangement satisfies $0 \leq c(i) < u_i$.

Note that the classical parking functions are recovered as ${\vec u}$-parking functions for the case ${\vec u} = (1,2, \dots, n)$. Also recall that classical parking functions coincide with the $G$-parking functions for the case $G = K_n = K_n^2$, the usual complete graph. For the case of complete hypergraphs we have the following generalization.

\begin{theorem}
\label{thm:ssuparking}
The $H$-parking functions of $K_{n+1}^d$ coincide with the set of ${\vec u}$-parking functions, where ${\vec u} = (u_1,\ldots,u_n)$ with:
\begin{itemize}
\item $ u_{k} =  \binom{n}{d-1} - \binom{n - k}{d-1}$, \; for $k = 1, 2, \dots, n+1-d$;
\item $ u_k = \binom{n}{d-1}$, \; for $k = n-d+2, n-d+3, \dots n$.
\end{itemize}
\end{theorem}

\begin{proof}
We let $[n]$ denote the non-sink vertices of $K_{n+1}^d$, and consider the degree of vertices $i \in T$ for various subsets $T \subset [n]$. On the one hand, if $|T| < d$ we have $\deg_T(i) = \binom{n}{d-1}$.  If $|T| \geq d$, then $\deg_T(i) = \binom{n}{d-1} - \binom{|T|-1}{d-1}$.

Now suppose $\vec{c}$ is an $H$-parking function for $K_{n+1}^d$, which we can assume is in weakly increasing order, so that $c_1 \leq c_2 \leq \cdots \leq c_n$. Suppose there exists an $k$ such that $c_k \geq u_k$, where $u_k$ is as described as above.  If we let $T = \{k, k+1, \dots, n\}$, we see that $\deg_T(i) \geq \deg_T(i)$ for all $i \in T$. This implies that $T$ is bounded by $\vec{c}$, a contradiction.

Now suppose $\vec{d}$ is a ${\vec u}$-parking function, where ${\vec u}$ is as above. The argument in the first paragraph shows that no set of vertices is bounded by $\vec{d}$, so that $\vec{d}$ is an $H$-parking function.
\end{proof}

Note that for any $\vec{u}$, the \emph{maximal} $\vec{u}$-vector parking functions are given by all permutations of the entries of $\vec{u}$. Hence by \cref{thm:ssuparking}, we see that the number of maximal $H$-parking functions for $K^d_{n+1}$ is given by $n!/(d-1)!$. From \cref{thm:maximalsup}, this value also counts the number of (certain) acyclic orientations of $K_{m,n}$, where $m = \binom{n+1}{d}$. The formula for the number of (all) acyclic orientations of an arbitrary complete bipartite graph is quite complicated \cite{Swenson}. As worked out by Cameron et al. in \cite{cameron2022}, this number is given by
\[a(K_{m,n}) = \sum_{j=1}^{\min\{m+1,n+1\}} (j-1)!^2 S(m+1,j) S(n+1, j),\]
where $S(i,j)$ denotes the Stirling number of the second kind.

A natural question to ask is if we have a way to determine the number of $H$-parking functions of a complete hypergraph.  From \cite{Yan15}, the number of ${\vec u}$-parking functions for any ${\vec u}$ can be determined by computing the determinant of a certain \emph{Steck matrix}. We recall this result next. 

\begin{theorem}\cite{StaPit}
The number PF($\vec{u}$) of ${\vec u}$-parking functions equals $n! \det D$, where $D$ is the $n \times n$ matrix with $ij$ entry
\[ \frac{u_i^{j-i+1}}{(j-i+1)!}\]
\noindent
if $j-i+1 \geq 0$ (and $0$ otherwise).
\end{theorem}

The ideas behind this formula go back to work of Steck \cite{Steck}.

\begin{example}
For $n=4$ and $d=3$, we have from \cref{thm:ssuparking} that $\vec{u} = (3,5,6,6)$, and hence the $H$-parking functions for $K_5^3$ are given by
\[\{(2,4,5,5), (2,4,5,4), \dots, (2,4,5,0), (2,4,4,5), (2,4,4,4),  \dots,(0,0,0,0)\}.\]
\noindent
The corresponding Steck matrix is given by
\renewcommand*{\arraystretch}{1.25}
\[D = \begin{pmatrix} 3 & \frac{9}{2} & \frac{27}{6} & \frac{81}{24} \\ 1 & 5 & \frac{25}{2} & \frac{125}{6} \\ 0 & 1 & 6 & \frac{36}{2} \\ 0 & 0 & 1 & 6 \end{pmatrix}. \]

From this we conclude that the number of $H$-parking functions of $K_5^3$ is
\[4! \det D = 1203.\]
\end{example}

This determinantal formula computes the number of $H$-parking functions for $K^d_{n+1}$, but it is not always easy to implement in practice. In \cite{Yan15} Yan has obtained closed formulas for PF($\vec{u}$) for certain values of $\vec{u}$, but unfortunately they do not apply here. Hence we ask the following.

\begin{question}
Is there a closed formula for the number of $H$-parking functions, where $H = K_n^d$ is the complete hypergraph?
\end{question}
   
It would also be interesting to find other hypergraphs $H$ with the property that the set of $H$-parking functions are recovered by ${\vec u}$-parking functions. A similar question for the case of $d=2$ (usual graphs) was considered by Gaydarov and Hopkins in \cite{GayHop}, where it was shown that this class is quite restrictive.

\section{Burning algorithms and tree like objects}
\label{sec:bijections}

In this section, we define our notion of spanning trees for hypergraphs, and provide a bijection between these objects and the set of $H$-parking functions. For this, we adapt the BCP algorithm from \cref{sec:burning} to the setting of hypergraphs by applying similar ideas to the underlying bipartite incidence graphs. It will be convenient to first translate some terminology to the setting of bipartite graphs.

Suppose $G$ is a bipartite graph, with vertex set parts $V(G) = L \sqcup R$.
In this case, the inverse BCP map (described in \cref{sec:burning}) that produces a spanning tree from a $G$-parking function $(b_1, \dots, b_n)$ can be understood in fewer steps. Instead of constructing the tree one vertex at a time, we can proceed in batches, moving from one side of the bipartite graph to the other.

In particular, at step $m$ of the process, instead of choosing a single vertex $p_m$ and edge $e_{p_m}$ to add to $t_{m-1}$ to create $t_m$, we choose all vertices in $v \in V_m$ that have the same height (in $t$) as $p_m$. We then add all such vertices and edges $e_v$ to create $t_{m'}$. We call this \emph{Bipartite Breadth-first Burning (BBB) algorithm}.  We then have the following observation.

\setlength{\columnsep}{0.5cm}
\begin{figure}[h]
\centering
\begin{multicols}{4}
    \begin{tikzpicture}[decoration={markings, mark= at position 0.5 with {\arrow{stealth}}}, scale = 1]
        \definecolor{color1}{RGB}{200,0,0}
        \definecolor{color2}{RGB}{0,0,200}
        \coordinate [label=right:\textcolor{black}{$q$}] (D) at (2,0);
        \coordinate [label=right:\textcolor{black}{$v_3$ \textcolor{color1}{0}}] (C) at (2,1);
        \coordinate [label=right:\textcolor{black}{$v_2$} \textcolor{color1}{1}] (B) at (2,2);
        \coordinate [label=right:\textcolor{black}{$v_1$} \textcolor{color1}{1}] (A) at (2,3);
        \coordinate [label=left:\textcolor{black}{$e_3$}] (F) at (0,0.5);
        \coordinate [label=left:\textcolor{black}{$e_2$}] (G) at (0,1.5);
        \coordinate [label=left:\textcolor{black}{$e_1$}] (H) at (0,2.5);
        \filldraw [black] (A) circle [radius=1.5pt]
                   (B) circle [radius=1.5pt]
                   (C) circle [radius=1.5pt]
                   (D) circle [radius=1.5pt]
                   (F) circle [radius=1.5pt]
                   (G) circle [radius=1.5pt]
                   (H) circle [radius=1.5pt];
        \draw [color2] (D) circle [radius=5pt];
        \draw [black] (H) to (A);
        \draw [black] (H) to (B);
        \draw [black] (H) to (C);
        \draw[color=black] (G) to (A)
                            (G) to (B)
                            (G) to (D);
        \draw[color=black] (F) to (A)
                            (F) to (C);
                           
        \draw[color=red] (F) to (D);
     \end{tikzpicture}
    \begin{tikzpicture}[decoration={markings, mark= at position 0.5 with {\arrow{stealth}}}, scale = 1]
        \definecolor{color1}{RGB}{200,0,0}
        \definecolor{color2}{RGB}{0,0,200}
        \coordinate [label=right:\textcolor{black}{$q$}] (D) at (2,0);
        \coordinate [label=right:\textcolor{black}{$v_3$ \textcolor{color1}{0}}] (C) at (2,1);
        \coordinate [label=right:\textcolor{black}{$v_2$} \textcolor{color1}{1}] (B) at (2,2);
        \coordinate [label=right:\textcolor{black}{$v_1$} \textcolor{color1}{1}] (A) at (2,3);
        \coordinate [label=left:\textcolor{black}{$e_3$}] (F) at (0,0.5);
        \coordinate [label=left:\textcolor{black}{$e_2$}] (G) at (0,1.5);
        \coordinate [label=left:\textcolor{black}{$e_1$}] (H) at (0,2.5);
        \filldraw [black] (A) circle [radius=1.5pt]
                   (B) circle [radius=1.5pt]
                   (C) circle [radius=1.5pt]
                   (D) circle [radius=1.5pt]
                   (F) circle [radius=1.5pt]
                   (G) circle [radius=1.5pt]
                   (H) circle [radius=1.5pt];
        \draw [color2] (D) circle [radius=5pt];
        \draw [black] (H) to (A);
        \draw [black] (H) to (B);
        \draw [black] (H) to (C);
        \draw[color=black] (G) to (A)
                            (G) to (B);
                            
        \draw[color=black] (F) to (A)
                            (F) to (C);
                           
        \draw[color=red] (F) to (D)
                        (G) to (D);
     \end{tikzpicture}
    \begin{tikzpicture}[decoration={markings, mark= at position 0.5 with {\arrow{stealth}}}, scale = 1]
        \definecolor{color1}{RGB}{200,0,0}
        \definecolor{color2}{RGB}{0,0,200}
        \coordinate [label=right:\textcolor{black}{$q$}] (D) at (2,0);
        \coordinate [label=right:\textcolor{black}{$v_3$ \textcolor{color1}{0}}] (C) at (2,1);
        \coordinate [label=right:\textcolor{black}{$v_2$} \textcolor{color1}{1}] (B) at (2,2);
        \coordinate [label=right:\textcolor{black}{$v_1$} \textcolor{color1}{1}] (A) at (2,3);
        \coordinate [label=left:\textcolor{black}{$e_3$}] (F) at (0,0.5);
        \coordinate [label=left:\textcolor{black}{$e_2$}] (G) at (0,1.5);
        \coordinate [label=left:\textcolor{black}{$e_1$}] (H) at (0,2.5);
        \filldraw [black] (A) circle [radius=1.5pt]
                   (B) circle [radius=1.5pt]
                   (C) circle [radius=1.5pt]
                   (D) circle [radius=1.5pt]
                   (F) circle [radius=1.5pt]
                   (G) circle [radius=1.5pt]
                   (H) circle [radius=1.5pt];
        \draw [color2] (D) circle [radius=5pt];
        \draw [black] (H) to (A);
        \draw [black] (H) to (B);
        \draw [black] (H) to (C);
        \draw[color=black] (G) to (A)
                            (G) to (B);
                            
        \draw[color=black] (F) to (A);

        \draw[color=red] (F) to (D)
                        (G) to (D)
                        (F) to (C);
     \end{tikzpicture}
    \begin{tikzpicture}[decoration={markings, mark= at position 0.5 with {\arrow{stealth}}}, scale = 1]
        \definecolor{color1}{RGB}{200,0,0}
        \definecolor{color2}{RGB}{0,0,200}
        \coordinate [label=right:\textcolor{black}{$q$}] (D) at (2,0);
        \coordinate [label=right:\textcolor{black}{$v_3$ \textcolor{color1}{0}}] (C) at (2,1);
        \coordinate [label=right:\textcolor{black}{$v_2$} \textcolor{color1}{1}] (B) at (2,2);
        \coordinate [label=right:\textcolor{black}{$v_1$} \textcolor{color1}{1}] (A) at (2,3);
        \coordinate [label=left:\textcolor{black}{$e_3$}] (F) at (0,0.5);
        \coordinate [label=left:\textcolor{black}{$e_2$}] (G) at (0,1.5);
        \coordinate [label=left:\textcolor{black}{$e_1$}] (H) at (0,2.5);
        \filldraw [black] (A) circle [radius=1.5pt]
                   (B) circle [radius=1.5pt]
                   (C) circle [radius=1.5pt]
                   (D) circle [radius=1.5pt]
                   (F) circle [radius=1.5pt]
                   (G) circle [radius=1.5pt]
                   (H) circle [radius=1.5pt];
        \draw [color2] (D) circle [radius=5pt];
        \draw [black] (H) to (A);
        \draw [black] (H) to (B);
        \draw [black] (H) to (C);
        \draw[color=black] 
                            (G) to (B);
                            
        \draw[color=black] (F) to (A);

        \draw[color=red] (F) to (D)
                        (G) to (D)
                        (F) to (C)
                        (G) to (A);
     \end{tikzpicture}
\end{multicols}
\caption{First four steps of the inverse BCP algorithm. Note that given the configuration $\vec{c} = (1,1,0)$, we build the corresponding spanning tree one edge at a time.}
\label{fig:BCPburning}

\setlength{\columnsep}{-7.75cm}
\begin{multicols}{2}
    \begin{tikzpicture}[decoration={markings, mark= at position 0.5 with {\arrow{stealth}}}, scale = 1]
        \definecolor{color1}{RGB}{200,0,0}
        \definecolor{color2}{RGB}{0,0,200}
        \coordinate [label=right:\textcolor{black}{$q$}] (D) at (2,0);
        \coordinate [label=right:\textcolor{black}{$v_3$ \textcolor{color1}{0}}] (C) at (2,1);
        \coordinate [label=right:\textcolor{black}{$v_2$} \textcolor{color1}{1}] (B) at (2,2);
        \coordinate [label=right:\textcolor{black}{$v_1$} \textcolor{color1}{1}] (A) at (2,3);
        \coordinate [label=left:\textcolor{black}{$e_3$}] (F) at (0,0.5);
        \coordinate [label=left:\textcolor{black}{$e_2$}] (G) at (0,1.5);
        \coordinate [label=left:\textcolor{black}{$e_1$}] (H) at (0,2.5);
        \filldraw [black] (A) circle [radius=1.5pt]
                   (B) circle [radius=1.5pt]
                   (C) circle [radius=1.5pt]
                   (D) circle [radius=1.5pt]
                   (F) circle [radius=1.5pt]
                   (G) circle [radius=1.5pt]
                   (H) circle [radius=1.5pt];
        \draw [color2] (D) circle [radius=5pt];
        \draw [black] (H) to (A);
        \draw [black] (H) to (B);
        \draw [black] (H) to (C);
        \draw[color=black] (G) to (A)
                            (G) to (B);
                            
        \draw[color=black] (F) to (A)
                            (F) to (C);
                           
        \draw[color=red] (F) to (D)
                        (G) to (D);
     \end{tikzpicture}
    \begin{tikzpicture}[decoration={markings, mark= at position 0.5 with {\arrow{stealth}}}, scale = 1]
        \definecolor{color1}{RGB}{200,0,0}
        \definecolor{color2}{RGB}{0,0,200}
        \coordinate [label=right:\textcolor{black}{$q$}] (D) at (2,0);
        \coordinate [label=right:\textcolor{black}{$v_3$ \textcolor{color1}{0}}] (C) at (2,1);
        \coordinate [label=right:\textcolor{black}{$v_2$} \textcolor{color1}{1}] (B) at (2,2);
        \coordinate [label=right:\textcolor{black}{$v_1$} \textcolor{color1}{1}] (A) at (2,3);
        \coordinate [label=left:\textcolor{black}{$e_3$}] (F) at (0,0.5);
        \coordinate [label=left:\textcolor{black}{$e_2$}] (G) at (0,1.5);
        \coordinate [label=left:\textcolor{black}{$e_1$}] (H) at (0,2.5);
        \filldraw [black] (A) circle [radius=1.5pt]
                   (B) circle [radius=1.5pt]
                   (C) circle [radius=1.5pt]
                   (D) circle [radius=1.5pt]
                   (F) circle [radius=1.5pt]
                   (G) circle [radius=1.5pt]
                   (H) circle [radius=1.5pt];
        \draw [color2] (D) circle [radius=5pt];
        \draw [black] (H) to (A);
        \draw [black] (H) to (B);
        \draw [black] (H) to (C);
        \draw[color=black] 
                            (G) to (B);
                            
        \draw[color=black] (F) to (A);

        \draw[color=red] (F) to (D)
                        (G) to (D)
                        (F) to (C)
                        (G) to (A);
     \end{tikzpicture}
\end{multicols}
\caption{First two steps of the BBB algorithm. This time we build the spanning tree in batches, adding multiple edges in each step.}
\label{fig:BBBburning}   
\end{figure}

\begin{proposition}\label{prop:BBBsameBCP}
Suppose $G$ is a bipartite graph where $V(G) = L \sqcup R$, and let $\vec{c}$ be a $G$-parking function. Then running the BBB algorithm and the inverse BCP algorithm output the same spanning tree $t$.
\end{proposition}

\begin{proof}
We show that starting at any point $t_{m-1}$ of the BCP algorithm and following the appropriate number of steps, we end up with the tree same as $t_m'$. 

Suppose that we are in step $m$ of the BCP algorithm, and let $V_m$, $t$, and $p_m \in V_m$ be as above.  Let $A = \{a_0 := p_m <_{\pi(t)} a_1 <_{\pi(t)} \dots <_{\pi(t)} a_b\}$ denote the vertices of $V_m$ having the same height (in $t$) as $p_m$ in $t$. Note that the set $A$ is contained in one part of the bipartite graph $G$.  Our first claim is that for all $i = 1, \dots, b$ we have $p_{m+i} = a_{m+i}$. To see this, note that the vertices of $V_{m+i} \setminus A$ for $i \geq 0$ have height strictly greater than $p_m$, and also that the breadth-first tree order refines the natural partial ordering coming from height within the tree. This establishes the claim.

Next we claim that the edge chosen in the $(m+i)$-th step of the BCP algorithm, which we denote $e_{p_{m+i}}$, is exactly same as the edge assigned to $a_{m+i} \in V_m$ in $t$ within the $m$-th step of the BCP algorithm, which we denote $e_{a_{m+i}}$. This follows from the observation that there are no edges among the vertices in $A$ (since they are contained in one part of $G$), and any edge connecting $a_{m+i}$ and $V_{m+i} \setminus \{a_0,\dots,a_i\}$ will be larger than $e_{a_{m+i}}$. The result follows.
\end{proof}

We illustrate the proof of \cref{prop:BBBsameBCP} in \cref{fig:BCPburning} and \cref{fig:BBBburning}. We see that the result of the first two steps of the inverse BCP algorithm is same as the first step of the BBB algorithm. Similarly, the first four steps of the inverse BCP algorithm is achieved by the first two steps of the BBB algorithm.

As discussed in \cite{ChebikinPylyavskyy05}, the BBB algorithm can be seen as a generalization of Dhar's algorithm (for bipartite symmetric digraphs) in the following sense.
Suppose $\vec{c} \in {\mathbb Z}_{\geq 0}^n$ is any vector. We mark the vertices of the graph $G$, starting with the root $q$. At each iteration of the algorithm, mark all vertices $v$ that have more marked neighbors than the value $\vec{c}_v$. If in the end all vertices are marked, then $\vec{c}$ is a $G$-parking function (equivalently, a superstable configuration on $G$). Conversely for every $G$-parking function, this algorithm marks all vertices. 

Given a $q$-rooted spanning tree $T$ of $G$, let $W_i$ be the set of vertices of $T$ of height $i$, where the height of a vertex is the distance to $q$ in $T$. We then have the following observation (compare to \cite[Proposition 5.1]{ChebikinPylyavskyy05}).

\begin{corollary}
\label{cor:wiheight}
$W_i$ is exactly the set of vertices marked at the $i$-th step of the BBB algorithm.
\end{corollary}

\subsection{Spanning trees}

We next define our notion of spanning trees for a hypergraph $H$. For this, we will work with a certain equivalence relation on (usual) spanning trees of bipartite graphs.

 \begin{definition}\label{def:equiv}
Suppose $B$ is a connected bipartite graph with vertex set $V(B) = L \sqcup R$, and with distinguished sink vertex $q \in R$. Spanning trees $T$ and $T^\prime$ of $B$ are said to be \emph{burning equivalent} if, when orienting the edges of $T$ and $T^\prime$ towards the sink, the set of edges directed from $R$ to $L$ is the same.
\end{definition}

See \cref{fig:burningtreelike} for an example of burning equivalent trees. With this we can define our notion of a spanning tree of a hypergraph (with choice of sink vertex $q$).

 \begin{definition}\label{def:tree}
Suppose $H$ is a hypergraph with sink vertex $q$, and let $B(H) = (\cV, \cE)$ denote its bipartite incidence graph. A \emph{($q$-rooted) spanning tree} of $H$ is a burning equivalence class $[T]$ of spanning trees of $B(H)$.
\end{definition}

\begin{remark}\label{rem:treeedges}
Note that the data of a spanning tree $[T]$ of $H$ is given by a subset $\LRT{T} \subset \cE$, namely the set of edges oriented $R$ to $L$ when orienting edges toward the sink $q$. This set $\LRT{T}$ has the property that every non-sink vertex $V(H) \setminus q$ is incident to some element of $\LRT{T}$. Also note that for any $T \in [T]$, the degrees of vertices on the left (from the set $E(H)$) will always be the same.
\end{remark}

\begin{figure}[h]
\hspace*{0cm}
   \centering
    \input{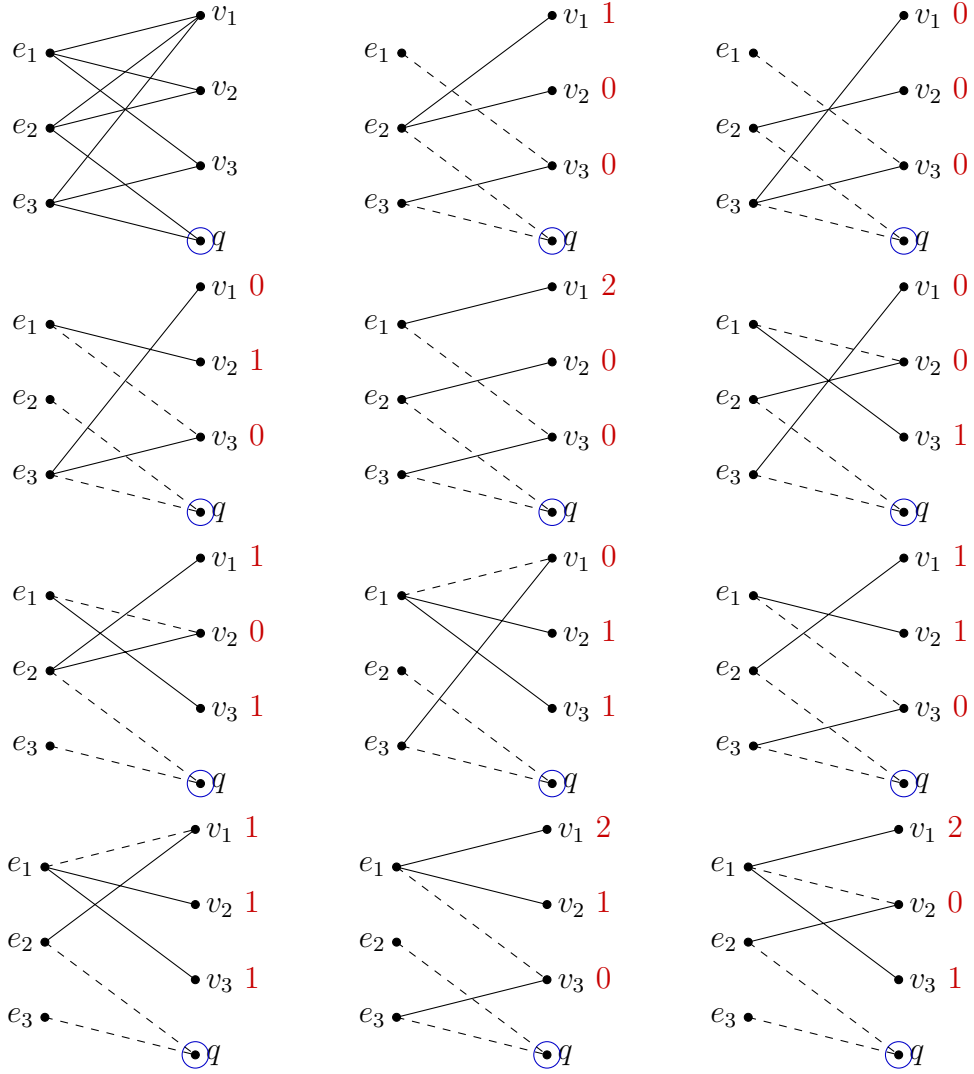}

    \caption{A hypergraph $H$ with representatives of each of its spanning trees, along with the corresponding $H$-parking functions. The solid edges represent the elements $\LRT{T}$ in each the burning equivalence class $[T]$, adding the dashed edges describe a spanning tree $T$ of $B(H)$ in this class.}
    \label{fig:treeclasses2}   
\end{figure}

\begin{remark}
We note that \cref{def:tree} generalizes the usual notion of a spanning tree of a graph. Indeed, if $G$ is a graph, then $B(G)$ has the property that each left-hand vertex $e \in E(G)$ has degree 2. A (classical) spanning tree $T$ of $G$ is then represented by a subgraph $B(T)$ of $B(G)$, where each vertex among $E(G)$ has degree 2, and where each vertex $v \in V(G)$ is incident to some edge in $B(T)$. To obtain an equivalence class of spanning trees for $B(G)$, we simply add the pendant degree 1 edges to the vertices among $E(G)$ that are not incident to any edge in $B(T)$. Note that any choice will lie in the same burning equivalence class since after orienting the edges toward $q$ these edges will be oriented left to right.

On the other hand, suppose $[T]$ is an equivalence class of spanning trees for $G$. Then for any representative $T$ (which is in particular a spanning tree of $B(G)$), the degrees among the vertices in $E(G)$ are preserved (see \cref{rem:treeedges}).
Let $\alpha$ and $\beta$ denote the number of vertices of degree 1 and degree 2 vertices of $T$ among the vertices of $E(G)$. Since $T$ is a spanning tree of $B(G)$, the total number of edges in $T$ is given by $|V(G)| + |E(G)| -1$ and hence $|V(G)| + |E(G)| -1 = \alpha + 2 \beta$.  But $\alpha = |E(G)| - \beta$ and hence $\beta = |V(G)| - 1$. We conclude that the set of degree $2$ vertices in $T \cap E(G)$ provides a well-defined spanning tree for $G$.
\end{remark}

\subsection{Bijections}

Now, suppose $B$ is a connected bipartite graph with vertex set ${\mathcal V} = L \sqcup R$ and specified sink vertex $q \in R$. 
Inspired by the breadth first tree orders described above, we will consider partitions of ${\mathcal V}$ that are refinements of $L \sqcup R$, so that $L = L_1 \sqcup L_2 \sqcup \cdots$ and $R = R_0 \sqcup R_1 \sqcup \cdots$.  
Given such a partition $W=(L_i,R_i)$, we use $W_{2i}$ to denote $R_i$ and $W_{2i-1}$ to denote $L_i$. We also set $W_0 = R_0 = \{q\}$.  We say that such a partition $W=(L_i,R_i)$ is \emph{valid} if 
\begin{enumerate}
    \item for each vertex $v$ in $W_i$, there is an edge in $B(H)$ that connects $v$ to some vertex in $W_{i-1}$, and
    \item for any vertex $w$ in $L_i$, there is no vertex $v \in R_j$ for $j = 0, 1 ,\dots i-2$ that is adjacent to $w$.
\end{enumerate}

We now apply these constructions to our setting. Suppose $H$ is a hypergraph with bipartite representation $B(H) = ({\mathcal V}, {\mathcal E})$, so that $L = E(H)$ and $R = V(H)$. Now, let $W=(L_i,R_i)$ be a valid partition of ${\mathcal V}$ as above. 
We let $T(W)$ denote the set of spanning trees of $B(H)$ that can be constructed by connecting each vertex in $W_i$ to some vertex in $W_{i-1}$. For any $T \in T(W)$, the vertices in $W_i$ have distance $i$ (in $T$) to the root vertex $q$. Hence, if the BBB algorithm applied to a given parking function on $B(H)$ results in this tree, the set $W_i$ is exactly the set of vertices that are marked at the $i$-th step. 

Note that for a fixed valid $W$, the height of a vertex $v$ has the same value for any $T \in T(W)$. Hence, for any $T \in T(W)$, we obtain the same ordering on the vertex set ${\mathcal V}$ coming from the tree order $\pi(T)$.

\begin{lemma}
\label{lem:treetolrdatainjective}
Suppose $\vec{c} \neq \vec{c}'$ are $H$-parking functions for $H$,
and let $T$ and $T^\prime$ be the corresponding spanning trees of $B(H)$ obtained by running the BBB algorithm on $\vec{c}_B$ and $\vec{c}'_B$. Then $T$ and $T^\prime$ are not burning equivalent (i.e. $\LRT{T} \neq \LRT{T'}$).
\end{lemma}

\begin{proof}
In the process of obtaining $T$ (resp. $T'$), we let $W_i$ (resp. $W^\prime_i$) be the set of vertices marked in the $i$th step of the BBB algorithm. From \cref{cor:wiheight}, we know that $W_i$ (resp. $W_i^\prime$) is also the set of vertices in $T$ (resp. $T^\prime$) of height $i$. Notice that from \cref{thm:CPtheorem}, $T$ and $T^\prime$ cannot be the same since $\vec{c} \neq \vec{c}^\prime$.

For the sake of contradiction, suppose that $\LRT{T} = \LRT{T^\prime}$. We will show by induction on $i$ that
\begin{itemize}
    \item $W_j = W^\prime_j$ for all $0 \leq j \leq i$, and
    \item $T|_{W_0 \cup W_1 \cup \cdots \cup W_i} = T^\prime|_{W^\prime_0 \cup W^\prime_1 \cup \cdots \cup W^\prime_i}$ for all $i$.
\end{itemize}
 For $i = 0$, this is clear. For $i$ even, we construct $W_{i+1} \subset L$ (resp. $W'_{i+1} \subset L$) as the set of all neighbors of $W_{i} \subset R$ (resp. $W'_i \subset R$) in the underlying graph $B(H)$ that have not yet appeared (recall that $\vec{c}_B$ and $\vec{c'}_B$ have all zeros assigned to the left vertices). Hence, we have $W_{i+1} = W'_{i+1}$.

Since $W_j = W'_j$ for $0 \leq j \leq i+1$, the tree order we get for $T$ and $T'$ assigns the same ordering to all edges connecting $W_{i+1}$ to $W_i$. Since we always pick the smallest possible edge in this case, we get $T|_{W_0 \cup W_1 \cup \cdots \cup W_{i+1}} = T^\prime|_{W^\prime_0 \cup W^\prime_1 \cup \cdots \cup W^\prime_{i+1}}$ as well.

For $i$ odd, from the definition of $\LRT{T}$, we have that $T|_{W_0 \cup W_1 \cup \cdots \cup W_{i+1}}$ is obtained from $T|_{W_0 \cup W_1 \cup \cdots \cup W_{i}}$ by adding the edges of $\LRT{T}$ that are incident to some element of $W_i$. By induction, we have $W_j = W'_j$ for $0 \leq j \leq i$, and by assumption we have $\LRT{T} = \LRT{T'}$. Hence, the claim follows.

The second part of the claim for tells us that $T = T^\prime$, which leads to a contradiction.
\end{proof}

We can now state the main result of this section.

\begin{theorem}\label{thm:bijection}
    For a hypergraph $H$ with root vertex $q$, there exists a bijection between the set of $H$-parking functions and the set of $q$-rooted spanning trees of $H$.
\end{theorem}
\begin{proof}
We use the BBB algorithm to define a function $\varphi$ from the set of $H$-parking functions of $H$ to the set of spanning trees of $H$. 
From \cref{lem:treetolrdatainjective} we see that $\varphi$ is injective. To show that $\varphi$ is surjective, let $T$ be a spanning tree of $B(H)$. 
It is enough to show that we can construct a spanning tree $T'$ of $B(H)$ such that $\LRT{T'} = \LRT{T}$, and that the $G$-parking function of $B(H)$ coming from $T'$ under the BBB algorithm is, in fact, an $H$-parking function for $H$ (that is, all left vertices are assigned zero chips).

For this, we first recursively construct a partition of ${\mathcal V}$ that is a refinement of $L \sqcup R$. We start with $W_0 := \{q\}$, and construct $W_1,W_2,W_3\dots$ as follows. For $k \geq 0$, we define $W_{2k+1}$ to be the set of neighbors of $W_{2k}$ that have not yet  been used in the partition. We define $W_{2k+2}$ to be the set of vertices in $R$ that have an edge in $\LRT{T'}$, and the opposite endpoint in $W_{2k+1}$. Observe that for any $k \geq 0$ we have $W_{2k+1} \subset L$, whereas $W_{2k+2} \subset R$. Then $W = W_0 \sqcup W_1 \sqcup \dots$ forms a valid partition of the set of vertices of $B(H)$. As discussed above, this defines a total ordering on ${\mathcal V}$ (and hence edges ${\mathcal E}$) by choosing the tree order defined by any element of $T(W)$.

For each vertex in $W_{2k+2}$, we consider the unique incident edge among $\LRT{T}$, using it to connect to a vertex of $L$. For each vertex in $W_{2k+1}$, we use the smallest edge (in the ordering described above) among the edges incident to a vertex in $W_{2k}$. This choice of edges produces a spanning tree $T'$. Note that by construction we have $\LRT{T'} = \LRT{T}$. 

Moreover, for each vertex $v \in W_{2k+1}$, there are no edges in $B(H)$ of the form $vw$ where $w \in W_0 \cup \dots \cup W_{2k-1}$ (otherwise $v$ will not be in $V_{2k+1}$ due to the way it is constructed). This means that when we use the BBB algorithm (in reverse) to get a corresponding parking function, the vertex $v$ is assigned $0$ chips. We conclude that $T'$ is the desired spanning tree of $B(H)$.
\end{proof}

\setlength{\columnsep}{-6cm}
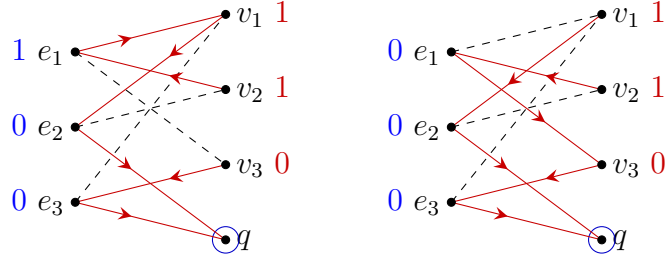
\begin{figure}[h]
\centering
\begin{multicols}{2}
    \begin{tikzpicture}[scale=1]
        \definecolor{color1}{RGB}{200,0,0}
        \definecolor{color2}{RGB}{0,0,200}
        \coordinate [label=right:\textcolor{black}{$q$}] (D) at (2,0);
        \coordinate [label=right:\textcolor{black}{$v_3$ \textcolor{color1}{0}}] (C) at (2,1);
        \coordinate [label=right:\textcolor{black}{$v_2$} \textcolor{color1}{1}] (B) at (2,2);
        \coordinate [label=right:\textcolor{black}{$v_1$} \textcolor{color1}{1}] (A) at (2,3);
        \coordinate [label=left:\textcolor{blue}{0} \textcolor{black}{$e_3$}] (F) at (0,0.5);
        \coordinate [label=left:\textcolor{blue}{0} \textcolor{black}{$e_2$}] (G) at (0,1.5);
        \coordinate [label=left:\textcolor{blue}{1} \textcolor{black}{$e_1$}] (H) at (0,2.5);
        \draw [color2] (D) circle [radius=5pt];
        \draw[color=black,dashed] (H) to (C);
        \draw[color=black,dashed] (G) to (B);
        \draw[color=black,dashed] (F) to (A);  
        \draw[color=color1, decoration={markings, mark= at position 0.375 with {\color{color1}{\arrow[line width=0.4mm]{stealth}}}}, postaction={decorate}] (F) to (D);
        \draw[color=color1, decoration={markings, mark= at position 0.375 with {\color{color1}{\arrow[line width=0.4mm]{stealth}}}}, postaction={decorate}] (C) to (F);
        \draw[color=color1, decoration={markings, mark= at position 0.375 with {\color{color1}{\arrow[line width=0.4mm]{stealth}}}}, postaction={decorate}] (H) to (A);
        \draw[color=color1, decoration={markings, mark= at position 0.375 with {\color{color1}{\arrow[line width=0.4mm]{stealth}}}}, postaction={decorate}] (B) to (H);
        \draw[color=color1, decoration={markings, mark= at position 0.375 with {\color{color1}{\arrow[line width=0.4mm]{stealth}}}}, postaction={decorate}] (A) to (G);
        \draw[color=color1, decoration={markings, mark= at position 0.375 with {\color{color1}{\arrow[line width=0.4mm]{stealth}}}}, postaction={decorate}] (G) to (D);
        \filldraw [black] (A) circle [radius=1.5pt]
                   (B) circle [radius=1.5pt]
                   (C) circle [radius=1.5pt]
                   (D) circle [radius=1.5pt]
                   (F) circle [radius=1.5pt]
                   (G) circle [radius=1.5pt]
                   (H) circle [radius=1.5pt];
     \end{tikzpicture}
    \columnbreak
    \begin{tikzpicture}[scale=1]
        \definecolor{color1}{RGB}{200,0,0}
        \definecolor{color2}{RGB}{0,0,200}
        \coordinate [label=right:\textcolor{black}{$q$}] (D) at (2,0);
        \coordinate [label=right:\textcolor{black}{$v_3$ \textcolor{color1}{0}}] (C) at (2,1);
        \coordinate [label=right:\textcolor{black}{$v_2$} \textcolor{color1}{1}] (B) at (2,2);
        \coordinate [label=right:\textcolor{black}{$v_1$} \textcolor{color1}{1}] (A) at (2,3);
        \coordinate [label=left:\textcolor{blue}{0} \textcolor{black}{$e_3$}] (F) at (0,0.5);
        \coordinate [label=left:\textcolor{blue}{0} \textcolor{black}{$e_2$}] (G) at (0,1.5);
        \coordinate [label=left:\textcolor{blue}{0} \textcolor{black}{$e_1$}] (H) at (0,2.5);
        \draw [color2] (D) circle [radius=5pt];
        \draw[color=black,dashed] (H) to (A);
        \draw[color=black,dashed] (G) to (B);
        \draw[color=black,dashed] (F) to (A);  
        \draw[color=color1, decoration={markings, mark= at position 0.375 with {\color{color1}{\arrow[line width=0.4mm]{stealth}}}}, postaction={decorate}] (F) to (D);
        \draw[color=color1, decoration={markings, mark= at position 0.375 with {\color{color1}{\arrow[line width=0.4mm]{stealth}}}}, postaction={decorate}] (C) to (F);
        \draw[color=color1, decoration={markings, mark= at position 0.625 with {\color{color1}{\arrow[line width=0.4mm]{stealth}}}}, postaction={decorate}] (H) to (C);
        \draw[color=color1, decoration={markings, mark= at position 0.375 with {\color{color1}{\arrow[line width=0.4mm]{stealth}}}}, postaction={decorate}] (B) to (H);
        \draw[color=color1, decoration={markings, mark= at position 0.625 with {\color{color1}{\arrow[line width=0.4mm]{stealth}}}}, postaction={decorate}] (A) to (G);
        \draw[color=color1, decoration={markings, mark= at position 0.375 with {\color{color1}{\arrow[line width=0.4mm]{stealth}}}}, postaction={decorate}] (G) to (D);
        \filldraw [black] (A) circle [radius=1.5pt]
                   (B) circle [radius=1.5pt]
                   (C) circle [radius=1.5pt]
                   (D) circle [radius=1.5pt]
                   (F) circle [radius=1.5pt]
                   (G) circle [radius=1.5pt]
                   (H) circle [radius=1.5pt];
     \end{tikzpicture}
\end{multicols}
\caption{Burning equivalent spanning trees $T$ and $T'$ of $B(H)$, along with the corresponding $G$-parking functions of $B(H)$ given by the BBB algorithm. Note that the configuration on the right corresponds to an $H$-parking function for the hypergraph $H$.}
\label{fig:burningtreelike}   
\end{figure}

\begin{example}
To illustrate the surjectivity of $\varphi$ in \cref{thm:bijection}, let $T$ denote the spanning tree of $B(H)$ depicted on the left side of \cref{fig:burningtreelike}.

Let $\beta$ be the total order on ${\mathcal V}$ (the vertices of $B(H)$) given by $q < v_3 < v_2 < v_1 < e_3 < e_2 < e_1$. We then have $W_0 = \{q\}, W_1 = \{e_2,e_3\}, W_2 = \{v_1,v_3\}, W_3 = \{e_1\}, W_4 = \{v_2\}$ as our valid partition coming from the height within $B(H)$.
Then, as discussed in the proof, the ordering $\beta$ and the partition $W$ give a total ordering on the vertices and edges of $B(H)$ from the tree ordering. For example, we have $e_1v_3 < e_1v_1$.

Recall that when constructing $T'$, we are take all the edges connecting $W_{2k+2}$ to $W_{2k+1}$ from $\LRT{T}$. Hence, in this case we add $(e_3,v_3),(e_2,v_1),(e_1,v_2)$ to $T'$. For the edges connecting $V_{2k+1}$ to $V_{2k}$, we take the smallest edges in the ordering mentioned above. For $e_3$ and $e_2$ we still take $e_3q$ and $e_2q$, respectively, but for $e_1$ we now choose $e_1v_3$, since $e_1v_3 < e_1v_1$. These choices result in the spanning tree $T'$ depicted on the right side of \cref{fig:burningtreelike}.
\end{example}

\subsection{Hypertrees}
In work of K\'alm\'an and Postnikov \cite{Kalman17}, the notion of a \emph{hypertree} associated to a hypergraph is introduced.  To recall this notion, suppose $H$ is a hypergraph with bipartite representation $B(H)$. A \emph{hypertree} in $H$ is a function (vector) $f: E \rightarrow {\mathbb N} = \{0,1, \dots\}$ with the property that there exists a spanning tree of $B(H)$ that has degree $f(e) + 1$ at each $e \in E$. In \cite{Kalman17} the authors develop a notion of \emph{internal activity} for these hypertrees that gives rise to a notion of a Tutte polynomial for polymatroids.

In our context, suppose $T$ is a spanning tree of $H$ (and hence, in particular, a subtree of $B(H)$). We obtain a hypertree $f:E \rightarrow {\mathbb N}$ from $T$ taking $f(e) = \deg_T(e) - 1$ if $\deg_T(e) > 0$ and $f(e) = 0$ if $\deg_T(e) = 0$. From \cref{rem:treeedges} we get the following observation.


\begin{proposition}\label{prop:hypertree}
For any hypergraph $H$ with sink $q$, two burning equivalent spanning trees give rise to the same hypertree.
\end{proposition}


For example, the spanning trees depicted in \cref{fig:burningtreelike} each correspond to the hypertree $(1,1,1)$. An important application of the bijection between spanning trees and superstable configurations preserving the external activity/degree. It would be interesting if we can extend the definition of external activity for spanning trees on hypergraphs and superstable configurations to our model so that it is preserved under the bijection described above.

\subsection{Star hypergraphs}
\label{sec:stars}

In this section we consider the special case of hypergraphs $H$ with the property that the sink $q$ vertex is in every edge of $H$, so that $q \in e$ for all $e \in E(H)$. We call these \emph{star hypergraphs}. It turns out that in this setting, the set of $H$-parking functions, as well as the bijection to spanning trees, can be understood purely in terms of the classical theory on a certain underlying digraph. This leads to a determinantal formula for the number of $H$-parking functions in this setting.

\begin{figure}
    \centering
    \begin{tikzpicture}[scale = 1]
        \definecolor{color1}{RGB}{200,0,0}
        \definecolor{color2}{RGB}{0,0,200}
        \coordinate [label=right:\textcolor{black}{$q$}] (D) at (2,0);
        \coordinate [label=right:\textcolor{black}{$v_3$}] (C) at (2,1);
        \coordinate [label=right:\textcolor{black}{$v_2$}] (B) at (2,2);
        \coordinate [label=right:\textcolor{black}{$v_1$}] (A) at (2,3);
        \coordinate [label=left:\textcolor{black}{$e_3$}] (F) at (0,0.5);
        \coordinate [label=left:\textcolor{black}{$e_2$}] (G) at (0,1.5);
        \coordinate [label=left:\textcolor{black}{$e_1$}] (H) at (0,2.5);
        \filldraw [black] (A) circle [radius=1.5pt]
                   (B) circle [radius=1.5pt]
                   (C) circle [radius=1.5pt]
                   (D) circle [radius=1.5pt]
                   (F) circle [radius=1.5pt]
                   (G) circle [radius=1.5pt]
                   (H) circle [radius=1.5pt];
        \draw [color2] (D) circle [radius=5pt];
        \draw [black, decoration={markings, mark= at position 0.375 with {\color{blue}{\arrow[line width=0.4mm]{stealth}}}}, postaction={decorate}] (A) to (H);
        \draw [black, decoration={markings, mark= at position 0.625 with {\color{blue}{\arrow[line width=0.4mm]{stealth}}}}, postaction={decorate}] (B) to (H);
        \draw [black, decoration={markings, mark= at position 0.625 with {\color{blue}{\arrow[line width=0.4mm]{stealth}}}}, postaction={decorate}] (C) to (H);
        \draw [black, decoration={markings, mark= at position 0.625 with {\color{red}{\arrow[line width=0.4mm]{stealth}}}}, postaction={decorate}] (H) to (D); 
        \draw [black, decoration={markings, mark= at position 0.625 with {\color{red}{\arrow[line width=0.4mm]{stealth}}}},postaction={decorate}] (G) to (D);
        \draw [black, decoration={markings, mark= at position 0.625 with {\color{red}{\arrow[line width=0.4mm]{stealth}}}},postaction={decorate}] (F) to (D);
        \draw[black, decoration={markings, mark= at position 0.375 with {\color{blue}{\arrow[line width=0.4mm]{stealth}}}}, postaction={decorate}] (A) to (G);
        \draw [black, decoration={markings, mark= at position 0.625 with {\color{blue}{\arrow[line width=0.4mm]{stealth}}}}, postaction={decorate}] (B) to (G);
        \draw [black,decoration={markings, mark= at position 0.4 with {\color{blue}{\arrow[line width=0.4mm]{stealth}}}}, postaction={decorate}] (A) to (F);
        \draw [black, decoration={markings, mark= at position 0.625 with {\color{blue}{\arrow[line width=0.4mm]{stealth}}}}, postaction={decorate}] (C) to (F);
     \end{tikzpicture}
    \caption{The digraph $\overleftarrow{H}$ associated to the star hypergraph $H$, where $H$ has edges $(123q, 12q, 13q)$. Note that $q$ appears in each edge of $H$.}
    \label{fig:directedhyper}
\end{figure}
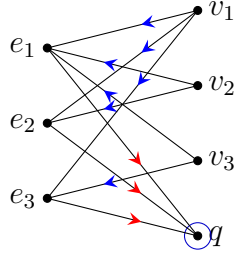

Suppose $H$ is a star hypergraph with sink vertex $q$ and nonsink vertices $[n] = \{1,2, \dots, n\}$. We define a digraph $\overleftarrow{H}$ by directing edges $(e,q)$ whenever $\{e,q\} \in \cE$ and $(v,e)$ whenever $\{v,e\} \in \cE$ and $v \neq q$. See \cref{fig:directedhyper} for an example. Note that in this example, $q \in e$ for all $e \in E(H)$ (so that this \emph{not} our running example).

\begin{lemma}
Suppose $H$ is a star hypergraph with sink vertex $q$, and with bipartite representation $B(H) = (\cV, \cE)$. Then there exists a bijection between the $q$-rooted spanning trees of $H$ (as a hypergraph) and the $q$-rooted spanning trees of $\overleftarrow{H}$ (as a digraph).
\end{lemma}

\begin{proof}
A $q$-rooted spanning tree of $\overleftarrow{H}$ is uniquely determined by a choice of output edge for each vertex. For vertices in $E(H)$, this is the unique edge incident to the the sink. For vertices in $V(H)$ (aside from $q$), this constitutes a unique choice of edge to a vertex connected to the sink. This is equivalent to choosing, for each $v_i \in V(H) \backslash q$, an edge oriented ``right to left" from $v_i$ to some element of $E(H)$. This is in turn equivalent to a choice of $q$-rooted spanning tree of $B(H)$, since the orientation on $\overleftarrow{H}$ is how we define burning equivalence on spanning trees of $H$.
\end{proof}

\begin{lemma}
Suppose $H$ is a hypergraph with sink vertex $q$, and with bipartite representation $B(H) = (\cV, \cE)$. Then the set of $G$-parking functions of $\overleftarrow{H}$ coincides with the set of $H$-parking functions on $H$.
\end{lemma}

\begin{proof}
Recall from \cref{prop:conversion} that the set of $H$-parking functions on $H$ coincide with the $G$-parking functions on $B(H)$ where all entries on coordinates from $E(H)$ are $0$. Now note that in $\overleftarrow{H}$, the degree of each vertex in $E(H)$ is $1$, implying that any $G$-parking function must have value $0$ there. On the other hand, the degree of any vertex among $V(H)$ coincides with its degree as a vertex in $B(H)$. The result follows.
\end{proof}

In the case of a star hypergraph, we also obtain a matrix (the reduced Laplacian of $\overleftarrow{H}$) whose determinant counts the number of $H$-parking functions . In our example from \cref{fig:directedhyper}, the matrix is given by

\[L_{\overleftarrow{H}} = \begin{pmatrix} 1 & 0 & 0 & 0 & 0 & 0 \\ 0 & 1 & 0 & 0 & 0 & 0 \\ 0 & 0 & 1 & 0 & 0 & 0 \\ -1 & -1 & -1 & 3 & 0 & 0 \\ -1 & -1 & 0 & 0 & 2 & 0 \\ -1 & 0 & -1 & 0 & 0 & 2 \end{pmatrix}\]
Here we have indexed rows and columns by $(e_1, e_2, e_3, v_1, v_2, v_3)$.

\begin{figure}[h]
    \centering
    \begin{subfigure}[b]{0.15\textwidth}
        \begin{tikzpicture}[decoration={markings, mark= at position 0.5 with {\arrow{stealth}}}, scale = 1]
        \definecolor{color1}{RGB}{200,0,0}
        \definecolor{color2}{RGB}{0,0,200}
        \coordinate [label=right:\textcolor{black}{$q$}] (D) at (2,0);
        \coordinate [label=right:\textcolor{black}{$v_3$}] (C) at (2,1);
        \coordinate [label=right:\textcolor{black}{$v_2$}] (B) at (2,2);
        \coordinate [label=right:\textcolor{black}{$v_1$}] (A) at (2,3);
        \coordinate [label=left:\textcolor{black}{$e_3$}] (F) at (0,0.5);
        \coordinate [label=left:\textcolor{black}{$e_2$}] (G) at (0,1.5);
        \coordinate [label=left:\textcolor{black}{$e_1$}] (H) at (0,2.5);
        \filldraw [black] (A) circle [radius=1.5pt]
                   (B) circle [radius=1.5pt]
                   (C) circle [radius=1.5pt]
                   (D) circle [radius=1.5pt]
                   (F) circle [radius=1.5pt]
                   (G) circle [radius=1.5pt]
                   (H) circle [radius=1.5pt];
        \draw [color2] (D) circle [radius=5pt];
        \draw [black] (H) to (A) (H) to (C) (H) to (B);
     \end{tikzpicture}
    \end{subfigure}
    \hspace{25mm}
    \begin{subfigure}[b]{0.15\textwidth}
        \begin{tikzpicture}[scale = 1]
        \definecolor{color1}{RGB}{200,0,0}
        \definecolor{color2}{RGB}{0,0,200}
        \coordinate [label=right:\textcolor{black}{$q$}] (D) at (2,0);
        \coordinate [label=right:\textcolor{black}{$v_3$}] (C) at (2,1);
        \coordinate [label=right:\textcolor{black}{$v_2$}] (B) at (2,2);
        \coordinate [label=right:\textcolor{black}{$v_1$}] (A) at (2,3);
        \coordinate [label=left:\textcolor{black}{$e_3$}] (F) at (0,0.5);
        \coordinate [label=left:\textcolor{black}{$e_2$}] (G) at (0,1.5);
        \coordinate [label=left:\textcolor{black}{$e_1$}] (H) at (0,2.5);
        \filldraw [black] (A) circle [radius=1.5pt]
                   (B) circle [radius=1.5pt]
                   (C) circle [radius=1.5pt]
                   (D) circle [radius=1.5pt]
                   (F) circle [radius=1.5pt]
                   (G) circle [radius=1.5pt]
                   (H) circle [radius=1.5pt];
        \draw [color2] (D) circle [radius=5pt];
        \draw [black, decoration={markings, mark= at position 0.65 with {\color{blue}{\arrow[line width=0.4mm]{stealth}}}}, postaction={decorate}] (A) to (H);
        \draw [black, decoration={markings, mark= at position 0.65 with {\color{blue}{\arrow[line width=0.4mm]{stealth}}}}, postaction={decorate}] (B) to (H);
        \draw [black, decoration={markings, mark= at position 0.65 with {\color{blue}{\arrow[line width=0.4mm]{stealth}}}}, postaction={decorate}] (C) to (H);
        \draw [black, decoration={markings, mark= at position 0.65 with {\color{red}{\arrow[line width=0.4mm]{stealth}}}},postaction={decorate}] (G) to (D);
        \draw [black, decoration={markings, mark= at position 0.65 with {\color{red}{\arrow[line width=0.4mm]{stealth}}}},postaction={decorate}] (F) to (D);
        \draw[black, decoration={markings, mark= at position 0.65 with {\color{red}{\arrow[line width=0.4mm]{stealth}}}}, postaction={decorate}] (H) to (D);
     \end{tikzpicture}
    \end{subfigure}
    \caption{A spanning tree of the hypergraph $H$, along with the associated spanning tree of $\overleftarrow{H}$.}
    \label{fig:missing tree}
\end{figure}
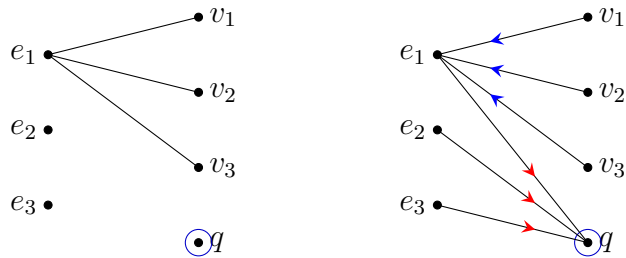

Note that $\det L_{\overleftarrow{H}} = 12$, and hence there are $12$ $q$-rroted spanning trees of $H$. Indeed, one can check that this set contains all 11 spanning trees depicted \cref{fig:treeclasses2}, along with the spanning tree in \cref{fig:missing tree}.

\begin{remark}
If $H$ is a star hypergraph with sink $q$ and nonsink vertices $[n]$, there is just one acyclic orientation of $H$ with unique source $q$, given by choosing $q \in e$ for all $e \in E(H)$. Hence by \cref{thm:maximalsup}, there exists a unique maximal $H$-parking function $\vec{c}$, given $\vec{c}_i = \deg(i)$. It follows that the number of $H$-parking functions is given by $\prod_{i=1}^n \deg(i)-1$.
\end{remark}

\begin{remark}
    Recall that a $q$-rooted spanning tree of a hypergraph $H$ is an equivalence class of spanning trees of the underlying bipartite graph $B(H)$. In \cref{thm:bijection} we construct a canonical representative for each class, with the property that under the BBB-algorithm this representative tree is sent to an $H$-parking function where all coordinates from $E(H)$ have value $0$. For the case of star hypergraphs, the representative trees are precisely the spanning trees of $B(H)$ where all vertices among $E(H)$ have height $1$. 
\end{remark}

\section{Chip-firing on hypergraphs}\label{sec:chipfiring}

We next seek to understand $H$-parking functions on a hypergraph $H = (V,E)$ in terms of a notion of chip-firing. For this we fix a sink vertex $q \in V$, let $[n]$ denote the set of non-sink vertices, and consider a configuration $\vec{c} \in {\mathbb Z}_{\geq 0}^n$ of chips. Recall that for a usual graph, when we fire a vertex $v \in [n]$, one chip is sent to each incident edge and then passed to the other vertex incident to that edge. For hypergraphs, we do something similar: When we fire a vertex $v$, it sends a chip to each of its incident hyperedges. But now within each of these hyperedges, we must make a \emph{choice} of which vertex receives the chip that comes in. We make this precise with the following definition.

\begin{definition}
    Suppose $H$ is a hypergraph and $v \in [n]$ is a vertex with incident edges $e_1, \ldots, e_k$. A \emph{firing choice} at $v$ consists of a choice of a vertex $v_i \in e_i$ where $v_i \neq v$, for each $i = 1, \dots, k$. We use $C_v = (e_1{\,:\,}v_1, \ldots, e_k{\,:\,}v_k)$ to denote a firing choice.
\end{definition}

Given a configuration $\vec{c}$, when we fire $v$, we think of sending a chip to each $e_i$, and then using the firing choice $C_v$ to pass that chip to $v_i$. Note that if $H$ is a (usual) graph, there is a canonical firing choice for each $v$ (namely, the other endpoint of each edge incident to $v$).

\begin{figure}[h]
\hspace*{-1.73 cm}
   \centering
    \begin{subfigure}[b]{0.15\textwidth}
            \centering
            \begin{tikzpicture}[scale = 1]
                \definecolor{color1}{RGB}{200,0,0}
                \definecolor{color2}{RGB}{25,161,21}
                \definecolor{color3}{RGB}{0,0,200}
                \coordinate (D) at (2,0);
                \coordinate [label=right:\textcolor{black}{$v_3$}] (C) at (2,1);
                \coordinate [label=right:\textcolor{black}{$v_2$}] (B) at (2,2);
                \coordinate [label=right:\textcolor{black}{$v_1$}] (A) at (2,3);
                \coordinate [label=right:\textcolor{black}{$v_4$}] (D') at (2.15,0);
                \coordinate [label=right:\textcolor{color1}{$1$}] (C') at (2.5,1);
                \coordinate [label=right:\textcolor{color1}{$0$}] (B') at (2.5,2);
                \coordinate [label=right:\textcolor{color1}{$1$}] (A') at (2.5,3);
                \coordinate [label=left:\textcolor{black}{$e_3$}] (E) at (0,0.5);
                \coordinate [label=left:\textcolor{black}{$e_2$}] (F) at (0,1.5);
                \coordinate [label=left:\textcolor{black}{$e_1$}] (G) at (0,2.5);
                \draw [color3,thick] (D) circle [radius=5pt];
                \draw [color1] (G) to (C)
                        (G) to (B) (G) to (A);
                \draw [color2] (F) to (B)
                        (F) to (D) (F) to (A);
                \draw [color3] (E) to (A)
                        (E) to (C) (E) to (D);
                \filldraw [black] (A) circle [radius=1.5pt]
                           (B) circle [radius=1.5pt]
                           (C) circle [radius=1.5pt]
                           (D) circle [radius=1.5pt]
                           (E) circle [radius=1.5pt]
                           (F) circle [radius=1.5pt]
                           (G) circle [radius=1.5pt];
            \coordinate [label=right:\textcolor{black}{$\xleftarrow{
             (\textcolor{color1}{e_1}: v_3, \;\textcolor{color2}{e_2}: v_4)
            }$}]  (H) at (2.85,1.5);
            \end{tikzpicture}
        \end{subfigure}\hspace{3.5cm}
    \begin{subfigure}[b]{0.15\textwidth}
            \centering
            \begin{tikzpicture}[scale = 1]
                \definecolor{color1}{RGB}{200,0,0}
                \definecolor{color2}{RGB}{25,161,21}
                \definecolor{color3}{RGB}{0,0,200}
                \coordinate (D) at (2,0);
                \coordinate [label=right:\textcolor{black}{$v_3$}] (C) at (2,1);
                \coordinate [label=right:\textcolor{black}{$v_2$}] (B) at (2,2);
                \coordinate [label=right:\textcolor{black}{$v_1$}] (A) at (2,3);
                \coordinate [label=right:\textcolor{black}{$v_4$}] (D') at (2.15,0);
                \coordinate [label=right:\textcolor{color1}{$0$}] (C') at (2.5,1);
                \coordinate [label=right:\textcolor{color1}{$2$}] (B') at (2.5,2);
                \coordinate [label=right:\textcolor{color1}{$1$}] (A') at (2.5,3);
                \coordinate [label=left:\textcolor{black}{$e_3$}] (E) at (0,0.5);
                \coordinate [label=left:\textcolor{black}{$e_2$}] (F) at (0,1.5);
                \coordinate [label=left:\textcolor{black}{$e_1$}] (G) at (0,2.5);
                \draw [color3,thick] (D) circle [radius=5pt];
                \draw [color1] (G) to (C)
                        (G) to (B) (G) to (A);
                \draw [color2] (F) to (B)
                        (F) to (D) (F) to (A);
                \draw [color3] (E) to (A)
                        (E) to (C) (E) to (D);
                \filldraw [black] (A) circle [radius=1.5pt]
                           (B) circle [radius=1.5pt]
                           (C) circle [radius=1.5pt]
                           (D) circle [radius=1.5pt]
                           (E) circle [radius=1.5pt]
                           (F) circle [radius=1.5pt]
                           (G) circle [radius=1.5pt];
            \coordinate [label=right:\textcolor{black}{$\xrightarrow{
                (\textcolor{color1}{e_1}: v_1, \; \textcolor{color2}{e_2}: v_1)
            }$}]  (H) at (2.85,1.5);
            \end{tikzpicture}
        \end{subfigure}\hspace{3.5cm}
    \begin{subfigure}[b]{0.15\textwidth}
            \centering
            \begin{tikzpicture}[scale = 1]
                \definecolor{color1}{RGB}{200,0,0}
                \definecolor{color2}{RGB}{25,161,21}
                \definecolor{color3}{RGB}{0,0,200}
                \coordinate (D) at (2,0);
                \coordinate [label=right:\textcolor{black}{$v_3$}] (C) at (2,1);
                \coordinate [label=right:\textcolor{black}{$v_2$}] (B) at (2,2);
                \coordinate [label=right:\textcolor{black}{$v_1$}] (A) at (2,3);
                \coordinate [label=right:\textcolor{black}{$v_4$}] (D') at (2.15,0);
                \coordinate [label=right:\textcolor{color1}{$0$}] (C') at (2.5,1);
                \coordinate [label=right:\textcolor{color1}{$0$}] (B') at (2.5,2);
                \coordinate [label=right:\textcolor{color1}{$3$}] (A') at (2.5,3);
                \coordinate [label=left:\textcolor{black}{$e_3$}] (E) at (0,0.5);
                \coordinate [label=left:\textcolor{black}{$e_2$}] (F) at (0,1.5);
                \coordinate [label=left:\textcolor{black}{$e_1$}] (G) at (0,2.5);
                \draw [color3,thick] (D) circle [radius=5pt];
                \draw [color1] (G) to (C)
                        (G) to (B) (G) to (A);
                \draw [color2] (F) to (B)
                        (F) to (D) (F) to (A);
                \draw [color3] (E) to (A)
                        (E) to (C) (E) to (D);
                \filldraw [black] (A) circle [radius=1.5pt]
                           (B) circle [radius=1.5pt]
                           (C) circle [radius=1.5pt]
                           (D) circle [radius=1.5pt]
                           (E) circle [radius=1.5pt]
                           (F) circle [radius=1.5pt]
                           (G) circle [radius=1.5pt];
            \end{tikzpicture}
        \end{subfigure}
    \caption{Example of different firing choices.}
    \label{fig:enter-label}
\end{figure}
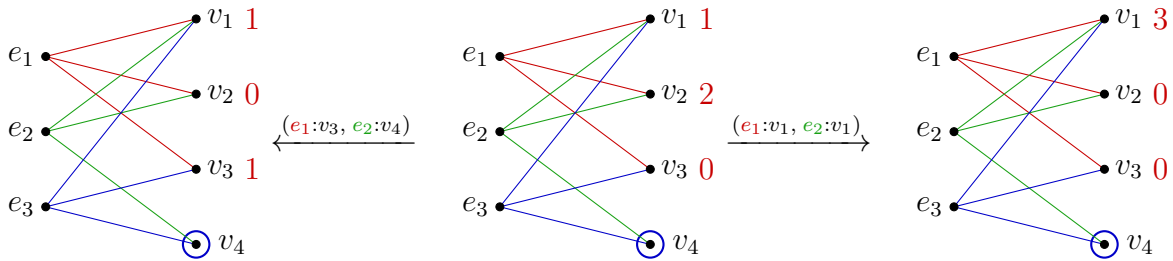

\begin{example}
Let $\vec{c} = (1,2,0)$ be the configuration on the hypergraph $H$ depicted in \cref{fig:enter-label}. We consider two different firing choices for the vertex $v_2$. If we fire $v_2$ with the firing choice $(e_1:v_3, e_2:v_4)$ we end up with the configuration $(1,0,1)$, whereas the choice $(e_1:v_1,e_2:v_1)$ gives us $(3,0,0)$.
\end{example}

We use $\cC_{v,H}$ to denote the set of all possible firing choices of the vertex $v$ in $H$. We write $\cC_v$ if the underlying hypergraph $H$ is clear. 
By considering all possible firing choices, we can recover the notion of when a vertex is ready to fire.

\begin{definition}
    Suppose $H$ is a hypergraph with sink vertex $q$ and nonsink vertices $[n]$, and let $\vec{c}$ be a configuration of chips. We say that a vertex $v \in [n]$ is \emph{ready to fire} if any firing choice $C \in \cC_v$ results in a nonnegative configuration.
\end{definition}

\begin{remark}
We note that any firing choice results in subtracting $\deg(v)$ from the coordinate of $\vec{c}$ corresponding to $v$ (and keeping all other entries at least as large as they were). Hence, a vertex is \emph{ready to fire} if it contains at least its degree in chips. It also follows that $v$ is ready to fire if \emph{some} firing choice results in a nonnegative configuration. 
\end{remark}

As was the case for usual graphs, we will also want to fire a set $T \subset V \setminus q$ of nonsink vertices. For this, we specify a firing choice for each vertex $v \in T$, but now we need an additional condition. Let $T \subset V \setminus q$ and suppose that $C_T = \{C_{v_i}: v_i \in T\}$ is a collection of firing choices for each element of $T$. Now suppose that there exists a hyperedge $e$ satisfying $e \subset T$.  
We say that $C_T$ is \emph{cancellative at $e$} if the firing choice $C_T$ restricted to $e$ results in a fixed point free permutation of the elements of $e$. Equivalently,
\[e = \{v_j\mid C_{v_i} = (\cdots, e:v_j, \cdots)\}_{v_i \in e}.\]
This generalizes set-firing in the classical graph case, where if two adjacent vertices are fired, a chip gets passed back and forth along that edge. Hence, the default firing choice for any chosen set to be fired is cancellative at every contained edge.

\begin{definition}
Suppose $H$ is a hypergraph, and let $T \subset V(H) \setminus q$. A firing choice $C_T$ is \emph{cancellative} (for $T$) if it is cancelative for all $e \subset T$.
\end{definition}

\begin{definition} \label{def:setfire}
    Suppose $H$ is a hypergraph and let $\vec{c}$ be a configuration on the non-sink vertices $V \setminus q$. A nonempty subset $T \subset V \setminus q$ is \emph{ready to fire} if every cancellative firing choice for $T$ results in a nonnegative configuration.
\end{definition}

\begin{example}\label{ex:choice}
For the hypergraph $H$ depicted in \cref{fig:enter-label}, consider the configuration $\vec{c} = (2,1,0)$, and the set $T = \{1,2,3\}$. An example of a cancellative firing choice is $\{C_{v_1} = (e_1:v_2, e_2:v_2, e_3:v_3), C_{v_2} = (e_1:v_3, e_2:v_4), C_{v_3} = (e_1:v_1, e_3:v_4)\}$, which results in a nonnegative configuration $(0,1,0)$. Another example of a cancellative firing choice is $\{C_{v_1} = (e_1:v_2, e_2:v_4, e_3:v_4), C_{v_2} = (e_1:v_3, e_2:v_4), C_{v_3} = (e_1:v_1, e_3:v_1)\}$, which results in the configuration $(1,0,-1)$. Hence $T$ is not ready to fire.  Note that both firing choices are indeed cancellative, since $e_1 \subset T$ and we have $e_1:v_1, e_1:v_2, e_1:v_3$ appearing exactly once each in each set.
\end{example}

With this we can define our notion of superstability.

\begin{definition}
    For a hypergraph $H$, a configuration $\vec{c}$ on the nonsink vertices is \emph{superstable} if no nonempty subset $T \subset V \setminus q$ is ready to fire.
\end{definition}

Although varying the cancellative firing choice of a set $T$ could lead to different distributions of chips, detecting whether $T$ is ready to fire can be detected in terms of degree. We record this as a lemma.

\begin{lemma}\label{lem:readytofire}
Suppose $\vec{c}$ is a configuration on $H$ and $T \subset V \setminus q$. Then $T$ is ready to fire if and only if $\deg_T(i) \leq c_i$ for all $i \in T$.
\end{lemma}

\begin{proof}
We consider the chip distribution on the vertex $i \in T$ when we fire the set $T$. Recall that for each edge $e$ incident to $i$, we send a chip to $e$. If $e \subset T$ then since the firing is cancellative, that edge will give a different chip back to $i$, resulting in no net change in chips. Hence, the total change in the number of chips at $i$ is at most $\deg_T(i)$. This shows that if $\deg_T(i) \leq c_i$ then $T$ is ready to fire. For the converse, if $\deg_T(i) > c_i$ for some $i \in T$, then we can come up with a cancellative firing choice that results in a configuration with a negative number of chips at $i$.
\end{proof}

Note that in \cref{ex:choice} we have $\deg_T(v_2) = 1$, so that $T$ was not ready to fire. On the other hand, if $\vec{c} = (2,1,1)$ then our condition implies that $T$ is ready to fire.

\begin{corollary}
Suppose $H$ is a hypergraph with sink vertex $q$.  Then a configuration $\vec{c}$ is superstable if and only if $\vec{c}$ is an $H$-parking function.
\end{corollary}

\begin{proof}
Recall that $\vec{c}$ is superstable if and only if no nonempty subset of the nonsink vertices is ready to fire. The result then follows from \cref{lem:readytofire} and the definition of an $H$-parking function.
\end{proof}

\subsection{Using digraphs to fire hypergraphs}\label{sec:digraphs}

As we have seen, firing a vertex $v$ in a hypergraph $H$ requires specifying a firing choice for each edge incident to $v$. One way to make this choice (for all vertices at once) is to cyclically order the elements of each hyperedge, giving rise to an underlying digraph on the vertices of $H$. In this section, we will show how chip-firing on this (collection of) digraphs can be used to recover superstable configurations on $H$ itself. In what follows, we will employ the theory of chip-firing on digraphs discussed in \cref{sec:digraphbackground}.  For this, we will need a way to construct a digraph from a hypergraph $H$. We start with a definition.

\begin{definition}
Let $H$ be a hypergraph. A \emph{cycling} $\cD$ of $H$ consists of a choice of cyclic orderings on the elements in each edge $e \in E(H)$.
\end{definition}

\begin{figure}[h]
\centering
\begin{multicols}{2}
    \vspace*{\fill}
    \begin{tikzpicture}[decoration={markings, mark= at position 0.5 with {\arrow{stealth}}}, scale = 1]
        \definecolor{color1}{RGB}{200,0,0}
        \definecolor{color2}{RGB}{25,161,21}
        \definecolor{color3}{RGB}{0,0,200}
        \definecolor{color4}{RGB}{0,0,0}
        \coordinate [label=right:\textcolor{black}{$v_4$}] (D) at (2,0);
        \coordinate [label=right:\textcolor{black}{$v_3$}] (C) at (2,1);
        \coordinate [label=right:\textcolor{black}{$v_2$}] (B) at (2,2);
        \coordinate [label=right:\textcolor{black}{$v_1$}] (A) at (2,3);
        \coordinate [label=left:\textcolor{black}{$e_3$}] (F) at (0,0.5);
        \coordinate [label=left:\textcolor{black}{$e_2$}] (G) at (0,1.5);
        \coordinate [label=left:\textcolor{black}{$e_1$}] (H) at (0,2.5);
        \draw[color=color1] (H) to (A)
                            (H) to (B)
                            (H) to (C);
        \draw[color=color2] (G) to (A)
                            (G) to (B)
                            (G) to (D);
        \draw[color=color3] (F) to (A)
                            (F) to (C)
                            (F) to (D);    
        \filldraw [black] (A) circle [radius=1.5pt]
                   (B) circle [radius=1.5pt]
                   (C) circle [radius=1.5pt]
                   (D) circle [radius=1.5pt]
                   (F) circle [radius=1.5pt]
                   (G) circle [radius=1.5pt]
                   (H) circle [radius=1.5pt];
     \end{tikzpicture}
       
        \label{fig:digraphfromhyper1}   
    \columnbreak
    \begin{tikzpicture}[decoration={markings, mark= at position 0.5 with {\arrow[line width=0.4mm]{stealth}}}, scale = 1.5]
        \definecolor{color1}{RGB}{200,0,0}
        \definecolor{color2}{RGB}{25,161,21}
        \definecolor{color3}{RGB}{0,0,200}
        \definecolor{color4}{RGB}{0,0,0}
        \coordinate [label=left:\textcolor{black}{$v_1$}] (A) at (0,1);
        \coordinate [label=above:\textcolor{black}{$v_2$}] (B) at (1,2);
        \coordinate [label=right:\textcolor{black}{$v_3$}] (C) at (2,1);
        \coordinate [label=below:\textcolor{black}{$v_4$}] (D) at (1,0);

        \draw[postaction={decorate}, color=color1] (B) to (A);
        \draw[postaction={decorate}, color=color2]  (A) to [bend left] (B) ;

        \draw[postaction={decorate}, color=color1] (C) to (B);

        \draw[postaction={decorate}, color=color3] (C) to (D);

        \draw[postaction={decorate}, color=color3] (D) to (A);
        \draw[postaction={decorate}, color=color2]  (D) to [bend left] (A) ;

        \draw[postaction={decorate}, color=color3] (A) to [bend right] (C);
        \draw[postaction={decorate}, color=color1]  (A) to [bend left] (C) ;

        \draw[postaction={decorate}, color=color2] (B) to [bend right] (D);

        \filldraw [black] (0,1) circle [radius=1.5pt]
                   (1,2) circle [radius=1.5pt]
                   (2,1) circle [radius=1.5pt]
                   (1,0) circle [radius=1.5pt];
        
\end{tikzpicture}
        
        \label{fig:digraphfromhyper2}
    \vfill

\end{multicols}

 \caption{A hypergraph $H$ and the digraph $D_\cD(H)$ associated to the cycling $\cD = \{(v_1v_3v_2), (v_1v_2v_4), (v_1v_3v_4)\}$.}
     \label{fig:digraphfromhyp}
\end{figure}
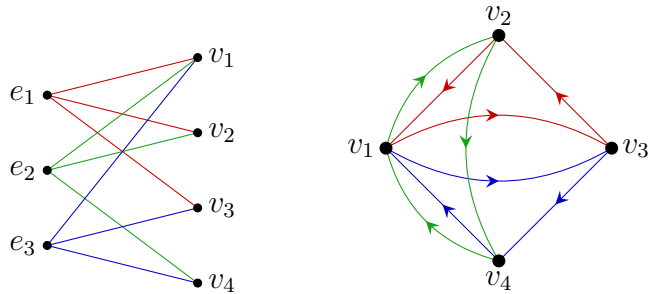

A cycling $\cD$ of $H$ defines a digraph $D_{\cD}(H)$ in the following way. The vertex set of $D_\cD(H)$ is taken to be $V(H)$, and we include the directed edge $(v_i, v_j)$ if $v_i \precdot_\cD v_j$ in some hyperedge $e$ of $H$. Here $x \precdot_\cD y$ indicates that $x$ immediately proceeds $y$ in the prescribed cyclic order on $e$.  See Figure \cref{fig:digraphfromhyp} for an example.

\begin{example}
For the hypergraph $H$ shown in \cref{fig:digraphfromhyp}, we consider the cycling
\[\cD = \{e_1 = (v_1 v_3 v_2), e_2 = (v_1v_2v_4), e_3 = (v_1v_3v_4)\}.\]
\noindent
The resulting induced digraph $D_\cD(H)$ is shown in \cref{fig:digraphfromhyp}. 
\end{example}

The digraphs we obtain from cyclings on $H$ lead to good notions of chip-firing, according to the following observation.

\begin{lemma}\label{lem:Eulerian}
    For any choice of cycling $\cD$ on a hypergraph $H$, the digraph $D_\cD(H)$ is Eulerian.
\end{lemma}

\begin{proof}
    Suppose $D_\cD(H)$ is a digraph obtained from a cycling $\cD$ of $H$. Let $v \in V$ be any vertex and consider the edges incident to $v$ in $D_\cD(H)$. Since each edge $e$ in $H$ is oriented cyclically, we have that if $v \in e$ then $e$ contributes one edge of the form $(v,v_i)$ and one of the form $(v_j,v)$ in $D_\cD(H)$ (where $v_i=v_j$ if and only if $|e| = 2$). Hence we have $\outdeg(v) = \indeg(v)$. To see that $D_\cD(H)$ is strongly connected, suppose $v$ and $w$ are vertices of $D_\cD(H)$. Since $H$ is assumed to be connected, we have a path $(v=v_0, v_1, \dots, v_k = w)$ through the hyperedges of $H$, where $\{v_i,v_{i+1}\} \subset e_i$ for some hyperedge $e_i \in E(H)$. Since each hyperedge $e \in E(H)$ defines a directed cycle in $D_\cD(H)$, we can construct a directed path from $v$ to $w$ by inserting a directed path from $v_i$ to $v_{i+1}$ (using a portion of the relevant cycle).
    \end{proof}

The collection of digraphs obtained from varying the choice of cyclings of $H$ provide a collection of superstable configurations. An important observation for us will be the following.

\begin{lemma}\label{lem:hyper to digraph}
Suppose $H$ is a hypergraph with nonsink vertices $[n]$ and suppose $\vec{c} \in {\mathbb Z}_{\geq 0}^n$ is a configuration of chips. If a nonempty subset $A \subset [n]$ is ready to fire in $H$, then for any cycling $\cD$ of $H$, $A$ is ready to fire in $D_\cD(H)$.
\end{lemma}

\begin{proof}
    If a subset $S \subset [n]$ is ready to fire in $H$, then by definition any firing choice leads to a a valid configuration. As $D(H)$ is one such way to send chips along the hyperedges, firing $A$ in $D(H)$ must result in a valid configuration as well. The result follows.
    \end{proof}

In what follows, we will consider cyclings $\cD$ on $H$ that are induced by an ordering of the vertex set $V(H)$. Note an ordering $v_1 < v_2 < \cdots < v_n$ of $V(H)$ induces a linear order on the elements of each edge $e_i = v_{i_1} < v_{i_2} < \dots < v_{i_p}$, which can be completed to the cycle $(v_{i_1} v_{i_2}  \dots v_{i_p})$. In this case, we say that $\cD$ is a \emph{vertex induced cycling}.

It turns out that one can obtain all superstable configurations for $H$ by considering vertex induced cyclings. More precisely we have the following result.

\begin{theorem}
\label{thm:hypersuperstable}
    Let $H$ be a hypergraph with fixed sink vertex $q$.
    A configuration $\vec{c}$ is superstable for $H$ if and only if it is superstable for $D_{\cD}(H)$ for some choice of vertex induced cycling $\cD$.
\end{theorem}

\begin{proof}
    For the forward direction, suppose that $\vec{c}$ is superstable for $H$. Then by definition the entire set of nonsink vertices $T = V\setminus q$ is not ready to fire. Hence there exists some cancellative firing choice $C_T$ that results in an invalid configuration. Without loss of generality, suppose $v_1$ is a vertex that ends up with a negative number chips as a result of this firing.
    
    We now consider firing the set of vertices $T_1 = V\setminus\{q,v_1\}$. Again since $\vec{c}$ is superstable, there exists some cancellative firing choice that results in a non-valid configuration. As before, let $v_2$ be a vertex that has a negative number of chips in the resulting invalid configuration. 
    We repeat this process until we have indexed all the vertices, and then use the resulting ordering $(q,v_1,v_2, \dots, v_n)$ to define a vertex induced cycling $\cD$.

    We now show that $\vec{c}$ is superstable as a configuration in the digraph $D_{\cD}(H)$ defined by $\cD$. It suffices to show that no nonempty subset of vertices is ready to fire. For this let $T\subset V\setminus q$ with $T \neq \emptyset$. Let $v_i \in T$ be the smallest element of $T$ under the ordering defined above. We then have that  $T\subset V^i$, where $V^i := V\setminus \{q, v_1,\dots, v_{i-1}\}$. 

We first show that $V^i$ is not ready to fire in $D_{\cD}(H)$.
    In the construction of our ordering, at the $i$-th step we have that set-firing $V^i$ in $H$ results in a negative number of chips at $v_i$. Recall that the firing of $V^i$ in $D_{\cD}(H)$ is one of the many firings of $V^i$ in $H$. We claim that the firing of $V^i$ in $D_{\cD}(H)$ results in a minimum number of chips at $v_i$ among all cancellative firings of $V^i$ in $H$. For this, we analyze how the chips are moving through each individual hyperedge containing $v_i$.

    To see this, first consider a hyperedge $e$ containing $v_i$ such that $e \subseteq V^i$. Since we only consider cancellative firings on $H$, any such such hyperedge $e$ should receive one chip from $v_i$ and also return one chip to $v_i$, resulting in $e$ not affecting the number of chips at $v_i$ before/after the firing. This is exactly what happens among the directed edges in $D_{\cD}(H)$ coming from $e$, when we consider the firing of $V^i$. 

    Next, suppose $e$ is a hyperedge containing $v_i$ with $e \nsubseteq V^i$. In $D_{\cD}(H)$, the only edge coming from $e$ that is directed toward $v_i$ is from a vertex $v \precdot_\cD v_i$. Since $e \nsubseteq V^i$, we have that $e$ must contain some vertex from the set $\{q,v_1,\dots,v_{i-1}\}$. We conclude that $v < v_i$. This means that $v \not \in V^i$, and hence $v$ is not fired and $v_i$ does not receive a chip from $e$. 
    
    The analysis of the two cases above implies that among all cancellative firings of $V^i$ in $H$, the firing of $V^i$ in $D_{\cD}(H)$ results in a minimum number of chips at $v_i$. In particular $v_i$ has a negative number of chips after firing $V^i$ in $D_{\cD}(H)$. 

    Now we go back to the arbitrary subset $T \subset V \backslash q$ and pick the maximum $i$ such that $T\subset V^i$. We claim that firing $T$ in $D_{\cD}(H)$ again results in a negative number of chips at $v_i$. To see this, note that $v_i \in T$ and hence $v_i$ still loses $\deg_{T}(v_i) \geq \deg_{V^i}(v_i)$ number of chips and it gains at most as many chips as it did when we fired $V^i$. Since $v_i$ had a negative number of chips when we fired $V^i$ in $D_{\cD}(H)$, the result follows.
    
    For the other direction we use \cref{lem:hyper to digraph}. In particular, if a subset cannot fire in $D_{\cD}(H)$, it cannot fire in $H$. Hence, the set of superstable configurations of $D_{\cD}(H)$ is always a subset of the superstable configurations of $H$.
\end{proof}

\begin{example}\label{ex:cycling}
Consider the hypergraph $H$ and digraph $D = D_\cD(H)$ depicted in \cref{fig:digraphfromhyp}, where the cycling is induced by the vertex ordering $v_1 < v_3 < v_2 < v_4$. The corresponding directed Laplacian is given by 
\[L_D = \begin{pmatrix}3 & -1 & -2 \\ -1 & 2 & 0 \\ 0 & -1 & 2 
\end{pmatrix}\]
One can check that the superstable configurations (which coincide with the set of $G$-parking functions) for $D$ are given by
\[\{000, 100,010, 001, 200, 110,101, 201\}.\]

This set provides a (proper) subset of the set of superstable configurations for the hypergraph $H$, which we have seen in \cref{fig:treeclasses2} to consist of
\[\{000,100,010,001,110,101,011,200,210,201,111\}.\]
\end{example}

\begin{remark}
\label{rmk:mincycles}
A natural question to ask is how many distinct cyclings on $H$ one must consider to produce all superstable configurations. Our \cref{thm:hypersuperstable} in particular implies that any superstable configuration $\vec{c}$ on $H$ can be obtained as a superstable configuration for the digraph $D_{\mathcal C}(H)$, for some choice of vertex induced cycling ${\mathcal C}$. Since the set of superstable configurations is closed under subtracting from a coordinate, this implies that the number of cyclings we must consider is at most the number of maximal superstable configurations on $H$. Recall that in \cref{sec:maximal} we provided a characterization of the maximal superstable configurations.
\end{remark}

\section{Further thoughts and open questions}

In this last section, we discuss some other aspects of parking functions on hypergraphs and also present some open questions. In this work we focused on combinatorial aspects of the theory, but there are also rich algebraic structures underlying these objects.

\subsection{Commutative algebra of hypergraph parking functions}
\label{sec:commalg}

As was the case for $G$-parking functions of classical graphs, one can define (and study) $H$-parking functions in an algebraic setting. Although these constructions are not needed for the results discussed above, they were in fact the original motivation for our study. We refer to \cite{AlmDocSmi} for details, where the underlying commutative algebra is the main focus.

\subsubsection{$H$-parking functions as standard monomials}
As was the case for usual graphs, the set of superstable configurations of a hypergraph $H$ is closed under the operation of subtracting `$1$' from a coordinate (assuming it does not become negative).  From this it follows that this set can be described as the \emph{standard monomials} of a certain monomial ideal. For graphs, this approach was adopted by Postnikov and Shapiro in \cite{PostnikovShapiro03}. A similar strategy can be employed in the hypergraph setting.

For a hypergraph $H$ with sink $q$ and non-sink vertices $[n] = \{1, \dots, n\}$, we fix a field ${\mathbb K}$ and let $S = {\mathbb K}[x_1, \dots, x_n]$ denote the polynomial ring.
If $T \subset [n]$ is a nonempty set of nonsink vertices, we let $m_T$ denote the monomial in $S$ given by
\[m_T = \prod_{i \in T} x_i^{\deg_T(i)},\]
\noindent
where $\deg_T(v)$ is as in \cref{def:outdegree}.

\begin{definition}
    Fix a field ${\mathbb K}$ and let $H$ be a hypergraph with sink $q$ and non-sink vertices $[n]$. The \emph{cut ideal} of $H$ is the monomial ideal in $S = {\mathbb K}[x_1, \dots, x_n]$ defined by
\[M_H = \langle m_T \rangle_{\emptyset \neq T \subset [n]}\]
\end{definition}

Recall that the \emph{standard monomials} of a monomial ideal $M$ are the monomials that do \emph{not} appear in $M$.  We then have the following observation.

\begin{lemma}
For a hypergraph $H$, the set of superstable configurations on $H$ corresponds to (exponent vectors of) the set of standard monomials of $M_H$.
\end{lemma}

\begin{proof}
    Suppose $\vec{c}$ is a configuration on the nonsink vertices of $H$, and let $T \subset [n]$ be a nonempty subset of the nonsink vertices. If we fire $T$ in a cancellative way, then each vertex $v \in T$ ends up losing $\deg_T(v)$ many chips. Hence after firing the set $T$ the resulting configuration has all of its coordinates nonnegative if and only if the monomial with exponent vector $\vec{c}$ is divisible by the monomial $m_T$.
\end{proof}

\begin{example}
If we take $H$ to be our running example depicted in \cref{fig:bipartite} (with sink $q = 4$), the ideal $M_H \subset {\mathbb K}[x_1, x_2, x_3]$ is given by
\[M_H = \langle x_1^3, x_2^2, x_3^2, x_1^2x_2x_3 \rangle.\]
\noindent
One can check that the exponent vectors of the associated standard monomials recover the superstable configurations of $H$. For instance, note that the maximal superstable configurations are given by $(2,1,0)$, $(2,0,1)$, $(1,1,1)$, which correspond to socle elements of $M_H$ that do not lie in $M_H$.
\end{example}

In \cite{PostnikovShapiro03}, Postnikov and Shapiro define $G$-parking functions (for a usual graph $G$) as the standard monomials of a similarly defined ideal $M_G$ (referred to as a \emph{tree ideal} in \cite{Klivans19}). They study homological properties of $M_G$, and for instance describe free resolutions and Betti numbers.

\subsubsection{Resolutions and generalized permutohedra}
\label{sec:genperm}
There is also a geometric approach to constructing the ideal $M_H$ which connects to the theory of generalized permutohedra.
This is related to (and inspired by) work of the second author and Sanyal from \cite{DocSan}, which we briefly recall here.

In \cite{DocSan}, the authors construct minimal cellular resolutions of the ideals $M_G$ for usual graphs $G$. The complex ${\mathcal B}_G$ supporting the resolution is obtained by intersecting an affine hyperplane $H_q$ (determined by choice of sink $q$) with the \emph{graphic arrangement} ${\mathcal A}_G$ associated to $G$. In particular, the minimal set of generators of $M_G$ are seen to coincide with the vertices (0-cells) of the complex ${\mathcal B}_G$.

If we dualize the above construction, we are lead to consider the \emph{graphic zonotope} ${\mathcal Z}_G$ associated to $G$. By definition ${\mathcal Z}_G$ is the Minkowski sum 
\[{\mathcal Z}_G = \sum_{ij \in E} [\vec{e}_i, \vec{e}_j]\]
\noindent
of line segments corresponding to edges of $G$. The normal fan of ${\mathcal Z}_G$ coincides with the graphic arrangement ${\mathcal A}_G$ (see \cite[Theorem 7.16]{Ziegler}).  Hence the vertices of ${\mathcal B}_G$ that one obtains by intersecting a hyperplane $H_q$ with ${\mathcal A}_G$ corresponds to the facets of ${\mathcal Z}_G$ that are visible from the direction normal to $H$.

In our context, given a hypergraph $H$ we construct the \emph{hypergraph polytope} ${\mathcal P}_H$, by definition the Minkowsi sum of simplices
\[{\mathcal P}_H = \sum_{e \in E(H)} \Delta_e,\]
\noindent
where $\Delta_e = \conv\{\vec{e}_i: i \in e\}$ is the simplex given by the standard basis vectors corresponding to elements of $e$. As before if we view ${\mathcal P}_H$ from the direction corresponding to $q$ we obtain a collection of facets that correspond to the generators of $M_H$. Even more, the complex supports a \emph{cocellular} resolutions of $M_H$, we refer to \cite{ADSRoot} and \cite{AlmDocSmi} for details.

\subsection{Open questions}
\label{sec:further}

We end with some open questions and possible directions for further research. 

Recall that for the case of a (usual) connected graph $G$, the degree sequence of $G$-parking functions can be interpreted in terms of activity, in the sense that the number of $G$-parking functions of degree $d$ is equal to the number of spanning trees with $d$ externally passive edges. In particular, work of Merino \cite{Merino01} shows that the degree sequence of $G$-parking functions recovers the $h$-vector of ${\mathcal M}(G)^*$, the matroid dual to the graphic matroid of $G$.

For the case of a hypergraph $H$, in \href{thm:bijection} we have shown that there exists a bijection between the set of $H$-parking functions and the set of spanning trees of $H$. However, we do not have any way to interpret the degree of an $H$-parking functions under this correspondence, and hence we ask the following question.

\begin{question}
    Is there some combinatorial interpretation for the degree sequence of $H$-parking functions on a hypergraph $H$? Is there some notion of activity for spanning trees of hypergraphs that can be defined for this?
\end{question}

For our next related question, results of Huh \cite{huh12} imply that the $h$-vector of a matroid defines a sequence that is \emph{log concave}. Hence for the case of graphs, the above discussion implies that the degree sequence of $G$-parking functions for any graph $G$ is log concave. At this point we do not know of any underlying matroid whose $h$-vector recovers the degree sequence of $H$-parking functions for hypergraphs. This inspires the following question.

\begin{question}
Is the degree sequence of $H$-parking functions of a hypergraph log concave?
\end{question}

Next, recall that in \cref{thm:hypersuperstable} we showed that the set of superstable configurations (equivalently the $H$-parking functions) of $H$ can be recovered by taking the union of all superstable configurations of the digraphs defined by all cyclings of $H$. This naturally leads to the following question, also discussed in \cref{rmk:mincycles} above.

\begin{question}
    Given a hypergraph $H$, can we describe a \emph{minimal} set of cyclings needed to recover all superstable configurations of $H$?
\end{question}

Associated to chip-firing on a graph $G$ is also the notion of a \emph{critical} configuration, by definition a configuration that is stable and \emph{reachable}.  The critical configurations can be used to define the \emph{critical group} of $G$.  In addition, the set of critical configurations is in a simple bijection with the set of superstable configurations. It is natural to ask if we have an analog of \cref{thm:hypersuperstable} for critical configurations.

\begin{question}
    Is there a natural notion of a critical configuration for a hypergraph $H$? If so can we find a bijection between these objects and the set of superstable configurations?
\end{question}

A natural candidate for a critical configuration on $H$ would be a critical configuration of the digraph $D_{\mathcal C}(H)$, where ${\mathcal C}$ is any cycling of $H$. We would then want an inherent definition of criticality for hypergraphs so that the set again is the union of all the critical configurations of digraphs associated to $H$. 
One subtlety regarding this question involves the fact that, in the classical setting, there are a number of characterizations of a critical configuration on a graph $G$ (see for instance \cite[Theorem 2.6.3]{Klivans19}). It turns out if we apply these various conditions on a hypergraph (or the associated bipartite graph) they yield different notions of criticality.

Our next questions address chip-firing on a hypergraph $H$. Recall that in our setup, when we fire a vertex $v$ one must make a choice of where each chip ends up.  Hence it is likely that any theory of chip-firing using these ideas must employ some stochastic or quantum aspects.

\begin{question}
Can one build a theory of chip-firing on a hypergraph $H$ using this model?
\end{question}

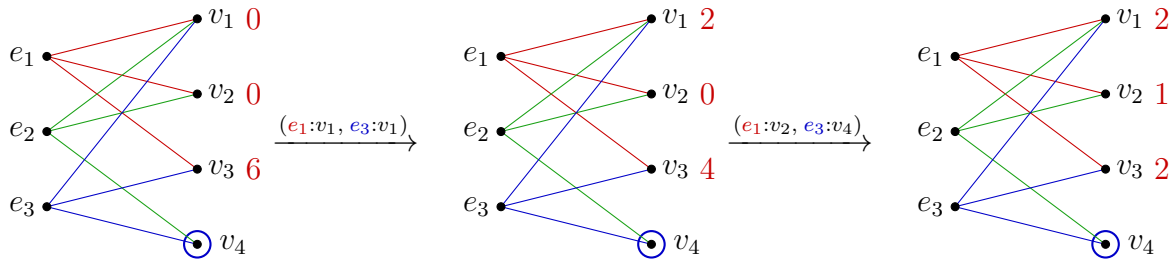
\begin{figure}[h]
\hspace*{-1.7 cm}
   \centering
    \begin{subfigure}[b]{0.15\textwidth}
            \centering
            \begin{tikzpicture}[scale = 1]
                \definecolor{color1}{RGB}{200,0,0}
                \definecolor{color2}{RGB}{25,161,21}
                \definecolor{color3}{RGB}{0,0,200}
                \coordinate (D) at (2,0);
                \coordinate [label=right:\textcolor{black}{$v_3$}] (C) at (2,1);
                \coordinate [label=right:\textcolor{black}{$v_2$}] (B) at (2,2);
                \coordinate [label=right:\textcolor{black}{$v_1$}] (A) at (2,3);
                \coordinate [label=right:\textcolor{black}{$v_4$}] (D') at (2.15,0);
                \coordinate [label=right:\textcolor{color1}{$6$}] (C') at (2.5,1);
                \coordinate [label=right:\textcolor{color1}{$0$}] (B') at (2.5,2);
                \coordinate [label=right:\textcolor{color1}{$0$}] (A') at (2.5,3);
                \coordinate [label=left:\textcolor{black}{$e_3$}] (E) at (0,0.5);
                \coordinate [label=left:\textcolor{black}{$e_2$}] (F) at (0,1.5);
                \coordinate [label=left:\textcolor{black}{$e_1$}] (G) at (0,2.5);
                \draw [color3,thick] (D) circle [radius=5pt];
                \draw [color1] (G) to (C)
                        (G) to (B) (G) to (A);
                \draw [color2] (F) to (B)
                        (F) to (D) (F) to (A);
                \draw [color3] (E) to (A)
                        (E) to (C) (E) to (D);
                \filldraw [black] (A) circle [radius=1.5pt]
                           (B) circle [radius=1.5pt]
                           (C) circle [radius=1.5pt]
                           (D) circle [radius=1.5pt]
                           (E) circle [radius=1.5pt]
                           (F) circle [radius=1.5pt]
                           (G) circle [radius=1.5pt];
            \coordinate [label=right:\textcolor{black}{$\xrightarrow{
             (\textcolor{color1}{e_1}: v_1, \; \textcolor{color3}{e_3}: v_1)
            }$}]  (H) at (2.85,1.5);
            \end{tikzpicture}
        \end{subfigure}\hspace{3.5cm}
    \begin{subfigure}[b]{0.15\textwidth}
            \centering
            \begin{tikzpicture}[scale = 1]
                \definecolor{color1}{RGB}{200,0,0}
                \definecolor{color2}{RGB}{25,161,21}
                \definecolor{color3}{RGB}{0,0,200}
                \coordinate (D) at (2,0);
                \coordinate [label=right:\textcolor{black}{$v_3$}] (C) at (2,1);
                \coordinate [label=right:\textcolor{black}{$v_2$}] (B) at (2,2);
                \coordinate [label=right:\textcolor{black}{$v_1$}] (A) at (2,3);
                \coordinate [label=right:\textcolor{black}{$v_4$}] (D') at (2.15,0);
                \coordinate [label=right:\textcolor{color1}{$4$}] (C') at (2.5,1);
                \coordinate [label=right:\textcolor{color1}{$0$}] (B') at (2.5,2);
                \coordinate [label=right:\textcolor{color1}{$2$}] (A') at (2.5,3);
                \coordinate [label=left:\textcolor{black}{$e_3$}] (E) at (0,0.5);
                \coordinate [label=left:\textcolor{black}{$e_2$}] (F) at (0,1.5);
                \coordinate [label=left:\textcolor{black}{$e_1$}] (G) at (0,2.5);
                \draw [color3,thick] (D) circle [radius=5pt];
                \draw [color1] (G) to (C)
                        (G) to (B) (G) to (A);
                \draw [color2] (F) to (B)
                        (F) to (D) (F) to (A);
                \draw [color3] (E) to (A)
                        (E) to (C) (E) to (D);
                \filldraw [black] (A) circle [radius=1.5pt]
                           (B) circle [radius=1.5pt]
                           (C) circle [radius=1.5pt]
                           (D) circle [radius=1.5pt]
                           (E) circle [radius=1.5pt]
                           (F) circle [radius=1.5pt]
                           (G) circle [radius=1.5pt];
            \coordinate [label=right:\textcolor{black}{$\xrightarrow{
                (\textcolor{color1}{e_1}: v_2, \;\textcolor{color3}{e_3}: v_4)
            }$}]  (H) at (2.85,1.5);
            \end{tikzpicture}
        \end{subfigure}\hspace{3.5cm}
    \begin{subfigure}[b]{0.15\textwidth}
            \centering
            \begin{tikzpicture}[scale = 1]
                \definecolor{color1}{RGB}{200,0,0}
                \definecolor{color2}{RGB}{25,161,21}
                \definecolor{color3}{RGB}{0,0,200}
                \coordinate (D) at (2,0);
                \coordinate [label=right:\textcolor{black}{$v_3$}] (C) at (2,1);
                \coordinate [label=right:\textcolor{black}{$v_2$}] (B) at (2,2);
                \coordinate [label=right:\textcolor{black}{$v_1$}] (A) at (2,3);
                \coordinate [label=right:\textcolor{black}{$v_4$}] (D') at (2.15,0);
                \coordinate [label=right:\textcolor{color1}{$2$}] (C') at (2.5,1);
                \coordinate [label=right:\textcolor{color1}{$1$}] (B') at (2.5,2);
                \coordinate [label=right:\textcolor{color1}{$2$}] (A') at (2.5,3);
                \coordinate [label=left:\textcolor{black}{$e_3$}] (E) at (0,0.5);
                \coordinate [label=left:\textcolor{black}{$e_2$}] (F) at (0,1.5);
                \coordinate [label=left:\textcolor{black}{$e_1$}] (G) at (0,2.5);
                \draw [color3,thick] (D) circle [radius=5pt];
                \draw [color1] (G) to (C)
                        (G) to (B) (G) to (A);
                \draw [color2] (F) to (B)
                        (F) to (D) (F) to (A);
                \draw [color3] (E) to (A)
                        (E) to (C) (E) to (D);
                \filldraw [black] (A) circle [radius=1.5pt]
                           (B) circle [radius=1.5pt]
                           (C) circle [radius=1.5pt]
                           (D) circle [radius=1.5pt]
                           (E) circle [radius=1.5pt]
                           (F) circle [radius=1.5pt]
                           (G) circle [radius=1.5pt];
            \end{tikzpicture}
        \end{subfigure}
    \caption{Example of different firing choices in a sequence of firing.}
    \label{fig:DiffFiringChoices}
\end{figure}

We emphasize that the firing choice at a vertex $v$ should be able to change each time we fire $v$ (since otherwise the theory reduces to chip-firing on a certain digraph). We illustrate this flexibility in firing choice in \cref{fig:DiffFiringChoices}.

\begin{example}
    Let $\vec{c} = (0,0,6)$ be the configuration on the hypergraph $H$ depicted on the left in \cref{fig:DiffFiringChoices}. 
    In the sequence of firings described there, we fire vertex 3 twice with different firing choice: first with $(123:1, 134:1)$ and then with $(123:2, 134:4)$.
\end{example}

Our last question addresses algebraic aspects of our study. Recall that for a hypergraph $H$, the $H$-parking functions can be recovered as the standard monomials of a certain monomial ideal $M_H$. In the classical setting of a graph $G$, the monomial ideal $M_G$ can be recovered as a certain initial ideal of the (binomial) toppling ideal $J_G$, by definition the lattice ideal defined by the reduced Laplacian of $G$. It is an open question weather such a structure exists in the context of hypergraphs.

\begin{question}
Is there a binomial \emph{toppling} ideal $J_H$ that encodes chip-firing on a hypergraph $H$ with the property that the monomial ideal $M_H$ arises as an initial ideal?
\end{question}

\section{Acknowledgements}
The work presented here was conducted as part of an REU at Texas State University in the summer of 2023, sponsored by NSF grant \#1757233. We thank NSF and Texas State for the support and the stimulating work environment. We are also grateful to Alex Constantinescu, Sophie Rehberg, and Ben Smith for helpful discussions.  Our definition of chip-firing on hypergraphs grew out of conversations that the second author had with Ben in early stages of the project. Dochtermann was partially supported by Simons Foundation Grant $\#964659$.

 \bibliographystyle{amsplain}
    \bibliography{references}

\end{document}